\documentclass[12pt,a4paper]{article}
\usepackage{amssymb}
\usepackage{graphicx}
\usepackage{amsmath}
\newtheorem{result}{Result}[section]

\newcounter{resultnum}[section]
\newtheorem{conclusion}{Conclusion}[section]

\newcounter{conclusionnum}[section]

\newcounter{conditionnum}[section]

\newcounter{conjecturenum}[section]
\newtheorem{example}{Example}[section]

\newcounter{examplenum}[section]

\newcounter{exercisenum}[section]

\newcounter{lemmanum}[section]

\newcounter{notationnum}[section]
\newtheorem{theorem}{Theorem}[section]

\newcounter{theoremnum}[section]
\newtheorem{definition}{Definition}[section]

\newcounter{definitionnum}[section]
\newtheorem{corollary}{Corollary}[section]

\newcounter{corollarynum}[section]
\newtheorem{remark}{Remark}[section]

\newcounter{remarknum}[section]
\newtheorem{proposition}{Proposition}[section]

\newcounter{propositionnum}[section]

\newcounter{acknowledgementnum}[section]

\newcounter{algorithmnum}[section]

\newcounter{axiomnum}[section]

\newcounter{casenum}[section]

\newcounter{claimnum}[section]

\newcounter{summarynum}[section]

\newcounter{problemnum}[section]
\newenvironment{proof}[1][]{\textbf{Proof.} }{}
\newcommand{ \R} {\mbox{\rm I$\!$R}}
\newcommand{ \C} {\mbox{\rm I$\!$C}}
\newcommand{ \Z} {\mbox{\rm I$\!$Z}}
\newcommand{ \I} {\mbox{\rm I$\!$I}}
\newcommand{ \D} {\mbox{\rm I$\!$\bf{D}}}

\begin{document}

\title{Nonholonomic  Clifford Structures 
   and Noncommutative Riemann--Finsler  Geometry}
\author{ Sergiu I. Vacaru\thanks{e--mail:
    vacaru@fisica.ist.utl.pt } \\ {---} \\  
Departamento de Fisica Aplicada,\\
 Facultad de Ciencias Experimentales,\\
 Universidad de Huelva, 21071,    Huleva, Spain \\
and \\
Departamento de  Matematica, Instituto Superior Tecnico, \\ 
 Av. Rovisco Pais1,  Lisboa, 1049--001, Portugal
   }
\maketitle

\begin{abstract}
We survey the geometry of Lagrange and Finsler spaces and discuss the
issues related to the definition of curvature of 
 nonholonomic manifolds enabled with  nonlinear
connection structure. It is proved that  any commutative
Riemannian geometry (in general, any Riemann--Cartan space) 
 defined by  a generic off--diagonal metric structure (with an
 additional  affine connection possessing  nontrivial torsion)
   is  equivalent  to a generalized
 Lagrange, or Finsler, geometry   modeled on 
 nonholonomic manifolds. 
This  results in  the problem of constructing
noncommutative  geometries with local anisotropy, 
in particular, related to geometrization of classical and quantum 
mechanical and field theories, even if we restrict our considerations only to 
commutative and noncommutative Riemannian spaces.  We  elaborate a 
geometric  approach  to the Clifford modules adapted to nonlinear connections,
 to the  theory of spinors and the Dirac operators  on nonholonomic
 spaces   and consider  possible  generalizations to 
 noncommutative geometry. We argue that any  commutative 
 Riemann--Finsler  geometry and  generalizations 
 my be derived from  noncommutative geometry by applying certain  
 methods elaborated for Riemannian spaces but extended to
 nonholonomic frame transforms and  manifolds provided with nonlinear
 connection structure.    

\vskip0.2cm 


\vskip0.2cm 

\end{abstract}

%\newpage

\tableofcontents

%\newpage

\section{ Introduction}

The goal of this work is to provide a better understanding of the
relationship between the theory of nonholonomic manifolds with
associated nonlinear connection structure, locally anisotropic spin
configurations and  Dirac operators on such manifolds and
noncommutative Riemann--Finsler and Lagrange geometry. The latter approach is
based on geometrical modelling of mechanical and classical field
theories (defined, for simplicity, by regular Lagrangians in analytic
mechanics and Finsler like anisotropic structures)
  and  gravitational, gauge and spinor field interactions in
low energy limits of string theory.  This allows to apply the
Serre--Swan theorem and think of vector bundles as projective modules,
which, for our purposes, are provided with nonlinear connection (in
brief, N--connection) structure and  can be defined as a
nonintegrabele (nonholonomic) distribution into conventional
horizontal  and vertical submodules. We relay on the theory of
 Clifford and spinor structures adapted to N--connections which
  results in locally anisotropic (Finsler like, or more general ones 
  defined by more general nonholonomic frame structures) Dirac
  operators. In the former item, it is the machinery of noncommutative
  geometry to derive distance formulas and to consider noncommutative
  extensions of Riemann--Finsler and Lagrange geometry and related off--diagonal
  metrics in gravity theories.

In \cite{vnc} it was proposed that an equivalent reformulation of the
general relativity theory as a gauge model with nonlinear realizations  
of the affine, Poincare and/or de Sitter groups allows a standard
extension of gravity theories in the language of noncommutative gauge
fields. The approach was developed in \cite{vncg} as an attempt to
generalize the A. Connes' noncommutative geometry \cite{connes1} to spaces with
generic local anisotropy. The nonlinear connection formalism was
elaborated for projective module spaces and the Dirac operator
associated to metrics in Finsler geometry and some generalizations
\cite{vsp1,vsph} (such as Sasaki type lifts of metrics to the tangent
bundles and vector bundle analogs) were considered as certain  
examples of noncommutative Finsler geometry. The constructions were
synthesized and revised in connection to ideas about appearance of
both noncommutative and Finsler geometry in string theory with
nonvanishing B--field and/or anholonomic (super) frame structures
 \cite{sw,con2,ard,abd,sah,vncgs,vstr0,vstr1} and in supergravity and gauge
 gravity  \cite{card,gars,vanc1,vanc2,dup}. In particular, one has 
 considered hidden noncommutative and Finsler like structures in
 general relativity and extra dimension gravity
 \cite{pan,vs1,v0,dv,v3}.

In this work, we confine ourselves to the classical aspects of
Lagrange--Finsler geometry (sprais, nonlinear connections, metric and
linear connection structures and almost complex structure derived from
from a Lagrange or Finsler fundamental form) in order  to 
generalize the  doctrine of the "spectral action" and the theory of
 distance in noncommutative geometry which is an extension of the 
previous results \cite{connes1}.  
For a complete information on modern noncommutative geometry and
physics,  we refer the reader to 
\cite{landi,madore,kastler2,bondia,douglas,konechny}, see a historical sketch in
Ref. \cite{kastler2} as well the aspects related to quantum group
theory \cite{manin,maj2,kassel} (here we note that the first quantum group Finsler
structure was considered in \cite{vargqugr}).  The theory of Dirac
operators and related Clifford analysis is a subject of various
investigations in modern mathematics and mathematical physics
\cite{martinmircea,martinmircea1,mitrea,nistor,fauser,castro1,castro2,castro3,castro4,schroeder,carey} 
(see also a relation to Finsler geometry \cite{vargcl} and an
off--diagonal "non" Kaluza--Klein compactified ansatz, but without
N--connection counstructions \cite{castro2a}). \footnote{The theory of
N--connections should not be confused with nonlinear gauge theories
and nonlinear relaizations of gauge groups.}  
 For an exposition spelling out all the
details of proofs and important concepts preliminary undertaken on 
  the subjects elaborated in our works, we refer to  proofs and
quotations in Refs. \cite{ma1,ma2,v1,v2,v3,vncg,vncgs,bondia,rennie1,lord,mart2}. 

This paper consists of two heterogeneous parts:

 The first (commutative) contains an overview of the Lagrange and Finsler
geometry and the off--diagonal metric and nonholonomic frame geometry in
gravity theories. In Section 2, we formulate the N--connection geometry
for arbitrary manifolds with tangent bundles admitting splitting into
conventional horizontal and vertical subspaces. We illustrate how regular
Lagrangians induce natural semispray, N--connection, metric and almost
complex structures on tangent bundles and discuss the relation between
Lagrange and Finsler geometry and theirs generalizations.
 We formulate six most important Results \ref{r1}--\ref{r6}
 demonstrating that the geometrization of Lagrange mechanics and the
 geometric models in gravity theories with generic off--diagonal metrics and
 nonholonomic frame structures are rigorously described by certain
 generalized Finsler geometries, which can be modeled equivalently
 both on Riemannian manifold and Riemann--Cartan nonholonomic
 manifolds. This give rise to the Conclusion \ref{ic} stating that a
 rigorous geometric study of nonholonmic frame and related metric,
 linear connection and spin structures in both commutative and
 noncommutative Riemann geometries requests the elaboration of
 noncommutative Lagrange--Finsler geometry.  Then, in
Section 3,  we consider  the theory of linear connections on N--anholonomic
manifolds (i. e. on manifolds with nonholonomic structure defined by
N--connections). We construct in explicit form the curvature
tensor of such spaces and define the Einstein equations for N--adapted
linear connection and metric structures.

The second (noncommutative)  part starts with Section 4 where we
define noncommutative N--anholonomic spaces. We consider the example
of noncommutative gauge theories adapted to the N--connection
structure.  Section 5 is devoted to
the geometry of nonholonomic Clifford--Lagrange structures. We define 
 the Clifford--Lagrange modules and Clifford N--anholonomic bundles
 being induced by the Lagrange quadratic form and adapted to the
 corresponding N--connection. Then we prove the {\bf Main Result 1,}
 of this work,   ( Theorem \ref{mr1}), stating that any regular
 Lagrangian and/or N--con\-nec\-ti\-on structure define naturally
 the fundamental geometric objects  and structures 
 (such as the Clifford--Lagrange module and Clifford d--modules, the 
 Lagrange spin structure  and d--spinors) for the corresponding Lagrange
 spin manifold  and/or  N--anholo\-nom\-ic spinor (d--spinor)
 manifold. We conclude that the Lagrange mechanics and
  off--diagonal gravitational interactions (in general, with nontrivial
  torsion and nonholonomic constraints) can be adequately 
  geometrized as certain Lagrange spin (N--anholonomic) manifolds. 

 In Section 6, we link up the theory of Dirac operators to
  nonholonomic structures
and spectral triples. We prove that there is a canonical spin
d--connection on the N--anholonomic manifolds generalizing that
induced by the Levi--Civita to the naturally ones induced by regular
Lagrangians and off--diagonal metrics. We define the Dirac d--operator
and the Dirac--Lagrange operator and formulate the {\bf Main Result 2}
(Theorem \ref{mr2}) arguing   that such N--adapted operators can be induced
canonically by almost Hermitian spin operators.  The concept of
distinguished spin triple is introduced in oder to  adapt the
constructions to the N--connection structure. Finally, the {\bf Main
  Result 3,} Theorem \ref{mr2}, is devoted to the definition, main
  properties  and  computation  of  distance in noncommutative spaces 
 defined by  N--anholonomic spin varieties. In these lecture notes, we
  only sketch in brief the ideas
  of proofs of the Main  Results:\ the details will be published in our
  further works.

\section{Lagrange--Finsler Geometry and Nonholonomic Manifolds}

This section presents some basic facts from the geometry of
nonholonomic manifolds  provided with nonlinear connection structure
 \cite{vsp1,nhm3,nhm1,nhm4,mironnh1,ver}.
 The constructions and methods are inspired from the Lagrange--Finsler
 geometry and generalizations
 \cite{fin,car1,rund,ma1,asa1,mats,bej,ma2,vstr,vstr1,bcs,mhss,sh} 
and gravity models on metric--affine 
spaces provided with generic off--diagonal metric, nonholonomic frame and affine
 connection structures \cite{vex1,vt,vs1,dv,v1,v2} 
(such spaces, in general, possess nontrivial torsion and
 nonmetricity).

\subsection{Preliminaries: Lagrange--Finsler metrics}\label{sslfm}

Let us consider a nondegenerate bilinear symmetric form $q(u,v)$ on a $n$--dimensional real vector space 
 $V^n.$ With respect to a basis $\{e_i\}^n_{i=1}$ for  $V^n,$ we express 
\begin{equation*}
q(u,v) \doteq q_{ij}u^i v^j
\end{equation*}%
for any vectors $u=u^i e_i,\ v= v^i e_i \in V^n$ and $q_{ij}$ being
  a nodegenerate symmetric matrix (the Einstein's
 convention on summing on repeating indices is adopted). This gives rise to 
 the Euclidean inner product 
 \begin{equation*} 
u \rfloor v \doteq q_{E}(u, v),
 \end{equation*}
if  $q_{ij}$ is positive definite, and  to the Euclidean norm  
 \begin{equation*} 
{|{\ \cdot \ } |}  \doteq \sqrt{q_{E}(u, u)}
 \end{equation*}
defining an Euclidean space $ (V^n, {|{\ \cdot \ } |}).$ Every
 Euclidean space
 is linearly isometric to the standard
 Euclidean space $\R ^n =  (R^n, {|{\ \cdot \ } |})$ if $q_{ij}=diag
 [1,1,...,1]$
 with standard Euclidean norm,
${|{\ y \ } |}  \doteq \sqrt{\sum_{i=1}^n |y^i|^2},$ for any
 $y=(y^i)\in R^n,$ where  $ R^n$ denotes the n--dimensional
 canonical real vector space.

There are also  different types of  quadratic forms/norms then the Euclidean
one:
\begin{definition} \label{dlagf}
A Lagrange  fundamental form $q_{L}(u, v)$\
on vector space  $V^n$ is defined by a
Lagrange  functional $L:\ V^n \to \R ,$ with  
\begin{equation}
q_{L(y)}(u,v)\doteq \frac{1}{2} \frac{\partial ^{2}}{\partial
  s\partial t}{[L(y+su+tu)]}_{\mid s=t=0} \label{lagrf}
\end{equation}
which is a ${\it C}^\infty$--function on $V^n \backslash \{0\}$ and 
 nondegenerate for any  nonzero vector $y \in V^n$ and  real
parameters $s$ and $t.$ 
\end{definition}
Having taken a  basis $\{e_i\}^n_{i=1}$ for  $V^n,$ we
transform $L=L(y^i e_i)$ in a function of $(y^i)\in R^n.$  

The Lagrange norm is ${|{\ \cdot \ } |}_L  \doteq \sqrt{q_{L}(u, u)}.$

\begin{definition} \label{dlargf} A Minkowski space is a pair $(V^n, F)$ where the
  Minkowski functional $F$ is a positively homogeneous of degree two
  Lagrange functional with the fundamental form (\ref{lagrf}) 
  defined for $L=F^2$ satisfying $F(\lambda y)= \lambda  F(y)$ for any
  $\lambda > 0$ and $y \in V^n.$
\end{definition}

The Minkowski norm is defined by  ${|{\ \cdot \ } |}_F  \doteq
\sqrt{q_{F}(u, u)}.$

\begin{definition}
The Lagrange (or Minkowski) metric fundamental function is defined 
\begin{equation}
g_{ij}=\frac{1}{2} \frac{\partial ^{2} L }{\partial
  y^i \partial y^j}(y) \label{lagm}
\end{equation}
(or \begin{equation}
g_{ij}=\frac{1}{2} \frac{\partial ^{2} F^2 }{\partial
  y^i \partial y^j} (y) \ ). \label{finm}
\end{equation}
\end{definition}

\begin{remark}\label{rem1} If $L$ is a Lagrange functional on $R^n$ (it could be
  also any functional of class  ${\it C}^\infty )$ with local
  coordinates  $(y^2,y^3,...,y^n),$ it also defines a singular
  Minkowski functional 
\begin{equation}
F(y)=\lbrack y^1 L(\frac{y^2}{y^1},...,\frac{y^n}{y^1}) \rbrack ^2
\end{equation}
which is of class  ${\it C}^\infty )$ on $R^n \backslash \lbrace y^1=0 
 \rbrace$.
\end{remark}

The Remark \ref{rem1} states that the Lagrange functionals are not
essentially more general than the Minkowski functionals
\cite{sh}. Nevertheless, we must introduce more general functionals
if we extend our considerations in relativistic optics, string models
of gravity and the theory of locally anisotropic stochastic and/or 
kinetic processes \cite{ma2,vex1,vt,vs1,dv,vkin}. 

Let us consider a base manifold $M,\ dim M =n,$ and its tangent bundle
$(TM, \pi , M)$ with natural surjective projection $\pi : TM \to M.$
From now on, all manifolds and geometric objects are supposed to be of
class $\mathcal{C}^{\infty }.$ We write $\widetilde{TM}= TM \backslash
 \lbrace 0 \rbrace $  where $ \lbrace 0 \rbrace $ means the null
 section of the map $\pi .$

A differentiable Lagrangian $L(x,y)$ is defined by a map
 $L:(x,y)\in TM\rightarrow L(x,y)\in \R$ of class $%
\mathcal{C}^{\infty }$ on $\widetilde{TM}$ and continous on the null section
$0:M\rightarrow TM$ of $\pi .$ For any point $x \in M,$ the
restriction $L_x \doteq L_{\mid T_xM}$ is a Lagrange functional on
$T_xM$ (see Definition \ref{dlagf}). For simplicity, in this work we
shall consider only  regular Lagrangians  with nondegenerated
Hessians,
\begin{equation}
\ ^{(L)}g_{ij}(x,y)=\frac{1}{2}\frac{\partial ^{2}L(x,y)}{\partial
y^{i}\partial y^{j}}  \label{lqf}
\end{equation}%
when $rank\left| g_{ij}\right| =n$ on $\widetilde{TM},$ which is a
Lagrange fundamental quadratic form (\ref{lagm}) on  $T_xM.$ In our further
considerations, we shall write $M_{(L)}$ if would be necessary to
emphasize that the manifold $M$ is provided in any its points with a
quadratic form (\ref{lqf}).

\begin{definition}
A Lagrange space is a pair $L^{n}=\left[ M,L(x,y)\right] $ with the
metric form  $\
^{(L)}g_{ij}(x,y)$ being of constant signature over $\widetilde{TM}.$
\end{definition}

\begin{definition}
A Finsler space is a pair $F^{n}=\left[ M,F(x,y)\right] $ where 
$F_{\mid x}(y)$ defines a Minkowski space with metric fundamental
function of type (\ref{finm}).
\end{definition}

The notion of Lagrange space was introduced by J. Kern \cite{kern} and
elaborated in details by the R. Miron's school on Finsler and Lagrange
geometry, see Refs. \cite{ma1,ma2}, as a natural extension of Finsler
geometry \cite{fin,car1,rund,asa1,mats,bej,bcs,sh} (see also Refs.\cite%
{vstr,vstr1}, on Lagrange--Finsler supergeometry, and
Refs. \cite{vnc,vncg,vncgs}, on some examples of noncommutative
locally anisotropic gravity and string theory).

\subsection{Nonlinear connection geometry}
We consider two important examples when the nonlinear connection (in
brief, N--connection) is
naturally defined in Lagrange mechanics and in gravity theories with
generic off--diagonal metrics and nonholonomic frames.

\subsubsection{Geometrization of mechanics: some important results}
\label{ssmr}
The Lagrange mechanics was geometrized by using methods of Finsler
geometry \cite{ma1,ma2} on tangent bundles enabled with a
corresponding nonholonomic structure (nonintegrable distribution)
 defining a canonical N--connection.\footnote{We cite a recent review
   \cite{ml} on alternative  approaches to geometric mechanics and
   geometry of classical fields related to investigation of  the
   geometric properties of
   Euler--Lagrange equations for various type of nonholonomic,
   singular or higher order systems. In the approach developed by
   R. Miron's school \cite{ma1,ma2,mhss}, the nonlinear connection and
 fundamental   geometric structures are derived in general form  from
 the  Lagrangian  and/or Hamiltonian: the basic  geometric
 constructions are not related to the  particular properties of
 certain  systems of partial  differential equations, 
 symmetries and constraints of mechanical  and field models.}
 By straightforward calculations, one proved the results:
\begin{result} \label{r1}
The Euler--Lagrange equations%
\begin{equation}
\frac{d}{d\tau }\left( \frac{\partial L}{\partial y^{i}}\right) -\frac{%
\partial L}{\partial x^{i}}=0 \label{eleq}
\end{equation}%
where $y^{i}=\frac{dx^{i}}{d\tau }$ for $x^{i}(\tau )$ depending on
parameter $\tau ,$ are equivalent to the ''nonlinear'' geodesic equations
\begin{equation}
\frac{d^{2}x^{i}}{d\tau ^{2}}+2G^{i}(x^{k},\frac{dx^{j}}{d\tau })=0
\label{ngeq}
\end{equation}%
where
\begin{equation}
2G^{i}(x,y)=\frac{1}{2}\ ^{(L)}g^{ij}\left( \frac{\partial ^{2}L}{\partial
y^{i}\partial x^{k}}y^{k}-\frac{\partial L}{\partial x^{i}}\right) 
\label{gcoeff}
\end{equation}%
with $^{(L)}g^{ij}$ being inverse to (\ref{lqf}).
\end{result}

\begin{result} \label{r2}
The coefficients $G^{i}(x,y)$ from (\ref{gcoeff}) define the solutions
of both type of equations (\ref{eleq}) and (\ref{ngeq})  as  paths of the canonical semispray%
\begin{equation*}
S=y^{i}\frac{\partial }{\partial x^{i}}-2G^{i}(x,y)\frac{\partial }{\partial
y^{i}}
\end{equation*}%
and a canonical N--connection structure  on $\widetilde{TM},$ 
\begin{equation}
^{(L)}N_{j}^{i}=\frac{\partial G^{i}(x,y)}{\partial y^{i}},  \label{cncl}
\end{equation}%
induced by the fundamental Lagrange function $L(x,y)$ (see Section
\ref{snc} on exact definitions and main properties).
\end{result}

\begin{result} \label{resws} 
The coefficients $\ ^{(L)}N_{j}^{i}$ defined by a Lagrange (Finsler)
fundamental function induce a global splitting on $TTM,$  a Whitney
sum, 
\begin{equation*}
TTM=hTM \oplus vTM 
\end{equation*}
as a nonintegrable distribution (nonholonomic, or equivalently,
anholonomic structure)into
horizontal (h) and vertical (v) subspaces parametrized
locally by frames (vielbeins) 
 $\mathbf{e}_{\nu }=(e_{i},e_{a}),$
where
\begin{equation}
e_{i}=\frac{\partial }{\partial x^{i}}-N_{i}^{a}(u)\frac{\partial }{\partial
y^{a}}\mbox{ and }e_{a}=\frac{\partial }{\partial y^{a}},  \label{dder}
\end{equation}%
and the dual frames (coframes) $\mathbf{\vartheta }^{\mu
}=(\vartheta ^{i},\vartheta ^{a}),$ where
\begin{equation}
\vartheta ^{i}=dx^{i}\mbox{ and }\vartheta ^{a}=dy^{a}+N_{i}^{a}(u)dx^{i}.
\label{ddif}
\end{equation}%
\end{result}

The vielbeins (\ref{dder}) and (\ref{ddif}) are called N--adapted (co)
frames. We omitted the label $(L)$ and used vertical indices
$a,b,c,...$  for the N--connection coefficients in order to be able to
use the formulas for arbitrary N--connections). We also note that we
shall use 'boldfaced' symbols for the geometric objects and spaces
adapted/ enabled to N--connection structure. For instance, we shall
write in brief $\mathbf{e}=(e,\ ^{\star }e)$ and
$\mathbf{\vartheta }=(\vartheta,\ ^{\star }\vartheta )$, respectively, for 
\begin{equation*}
\mathbf{e}_{\nu }=(e_{i},\ ^{\star }e_{k})=(e_{i},\  e_{a})
\mbox{ and } 
\mathbf{\vartheta }^{\mu}=(\vartheta ^{i},\ ^{\star }\vartheta ^{k})=
(\vartheta ^{i}, \vartheta ^{a}). 
\end{equation*}
The vielbeins (\ref{dder}) satisfy the nonholonomy relations
\begin{equation}
\lbrack \mathbf{e}_{\alpha },\mathbf{e}_{\beta }]=\mathbf{e}_{\alpha }%
\mathbf{e}_{\beta }-\mathbf{e}_{\beta }\mathbf{e}_{\alpha }=W_{\alpha \beta
}^{\gamma }\mathbf{e}_{\gamma }  \label{anhrel}
\end{equation}%
with (antisymmetric) nontrivial anholonomy coefficients $W_{ia}^{b}=\partial
_{a}N_{i}^{b}$ and $W_{ji}^{a}=\Omega _{ij}^{a}$ where 
\begin{equation}
\Omega _{ij}^{a}=\delta _{\lbrack j}N_{i]}^{a}=\frac{\partial N_{i}^{a}}{%
\partial x^{j}}-\frac{\partial N_{j}^{a}}{\partial x^{i}}+N_{i}^{b}\frac{%
\partial N_{j}^{a}}{\partial y^{b}}-N_{j}^{b}\frac{\partial N_{i}^{a}}{%
\partial y^{b}}.  \label{ncurv}
\end{equation}

 In order to preserve a relation with our previous denotations \cite%
{vex1,vsp1,vstr1}, we note that $\mathbf{e}_{\nu
}=(e_{i},e_{a})$ and $\mathbf{\vartheta }^{\mu }=(\vartheta ^{i},\vartheta
^{a})$ are, respectively, the former $\delta _{\nu }=\delta /\partial u^{\nu
}=(\delta _{i},\partial _{a})$ and $\delta ^{\mu }=\delta u^{\mu
}=(dx^{i},\delta y^{a})$ which emphasize that the  operators (\ref{dder}) and (%
\ref{ddif}) define, correspondingly, certain 'N--elongated' partial
derivatives and differentials which are more convenient for calculations on
spaces provided with nonholonomic structure. 

\begin{result} \label{r4}
 On $\widetilde{TM},$ there is a canonical metric structure 
$\ ^{(L)}\mathbf{g}= [g,\ ^{\star}g],$  
\begin{equation}
\ ^{(L)}\mathbf{g}=\ ^{(L)}g_{ij}(x,y)\ \vartheta ^{i}\otimes \vartheta
^{j}+\ ^{(L)}g_{ij}(x,y)\ ^{\star }\vartheta ^{i}\otimes \ ^{\star }\vartheta
^{j}  \label{slm}
\end{equation}%
constructed as a Sasaki type lift from $M.$\footnote{In
  Refs. \cite{sh,mhss}, it was suggested to use lifts with h- and
  v--components of type
 $\ ^{(L)}\mathbf{g}=(g_{ij}, g_{ij}a/\parallel y \parallel)$
  where $a=const$ and $\parallel y \parallel = g_{ij}y^i y^j $ in order
  to elaborate more physical extensions of the general relativity to
  the tangent bundles of manifolds. In another turn, such
  modifications are not necessary if we model Lagrange--Finsler
  structures by exact solutions with generic off--diagonal metrics in
  Einstein and/or gravity \cite{vex1,vt,vs1,dv,v2,v3,vkin}. For
  simplicity, in this work, we  consider only lifts of metrics of
  type (\ref{slm}).} 
\end{result}
We note that a complete geometrical model of Lagrange mechanics or a
well defined Finsler geometry can be elaborated only by additional assumptions
about a linear connection structure, which can be adapted, or not, to
a defined N--connection (see Section \ref{sdlc}).

\begin{result} \label{r5}
 The canonical N--connection (\ref{cncl})
induces naturally an almost complex structure $\mathbf{F}:\chi (\widetilde{TM%
})\rightarrow \chi (\widetilde{TM}),$ where $\chi $ denotes the module of
vector fields on $\widetilde{TM},$%
\begin{equation*}
\mathbf{F}(e_{i})=\ ^{\star }e_{i}\mbox{ and }\mathbf{F}(\ ^{\star
}e_{i})=-e_{i},
\end{equation*}%
when
\begin{equation}
\mathbf{F}=\ ^{\star }e_{i}\otimes \vartheta ^{i}-e_{i}\otimes \ ^{\star
}\vartheta ^{i}  \label{acs1}
\end{equation}%
satisfies the condition $\mathbf{F\rfloor \ F=-I,}$ i. e. $F_{\ \ \beta
}^{\alpha }F_{\ \ \gamma }^{\beta }=-\delta _{\gamma }^{\alpha },$ where $%
\delta _{\gamma }^{\alpha }$ is the Kronecker symbol and ''$\mathbf{\rfloor }
$'' denotes the interior product.
\end{result}
The last result is important for elaborating  an approach to geometric
quantization of mechanical systems modeled on nonholonomic manifolds
\cite{esv}  as well for definition of almost complex structures
derived from the real N--connection geometry related to nonholonomic
(anisotropic)  Clifford structures and  spinors in commutative
 \cite{vsp1,vsph,vst,vsp2,vp} and noncommutative spaces \cite{vnc,vncg,vncgs}.

\subsubsection{N--connections in  gravity theories}

For nonholonomic geometric models of gravity and string theories, one  does not consider  the
bundle $\widetilde{TM}$ but a general manifold $\mathbf{V},\ 
dim\mathbf{V}=n+m,$ which is a (pseudo) Riemannian space or a certain
generalization with possible torsion and nonmetricity fields.  
 A metric structure is defined on $\mathbf{V},$ with the coefficients
 stated with respect to a local coordinate basis $du^{\alpha }=\left(
dx^{i},dy^{a}\right) ,$ \footnote{the indices run correspondingly the values
 $i,j,k,...=1,2,...,n$ and $a,b,c,...=n+1,n+2,...,n+m.$}
\begin{equation*}
\mathbf{g}=\underline{g}_{\alpha \beta }(u)du^{\alpha }\otimes du^{\beta }
\end{equation*}%
where
\begin{equation}
\underline{g}_{\alpha \beta }=\left[
\begin{array}{cc}
g_{ij}+N_{i}^{a}N_{j}^{b}h_{ab} & N_{j}^{e}h_{ae} \\
N_{i}^{e}h_{be} & h_{ab}%
\end{array}%
\right] .  \label{ansatz}
\end{equation}

A metric, for instance, parametrized in the form (\ref{ansatz}), is generic
off--diagonal if it can not be diagonalized by any coordinate transforms.
Performing a frame transform with the coefficients
\begin{eqnarray}
\mathbf{e}_{\alpha }^{\ \underline{\alpha }}(u) &=&\left[
\begin{array}{cc}
e_{i}^{\ \underline{i}}(u) & N_{i}^{b}(u)e_{b}^{\ \underline{a}}(u) \\
0 & e_{a}^{\ \underline{a}}(u)%
\end{array}%
\right] ,  \label{vt1} \\
\mathbf{e}_{\ \underline{\beta }}^{\beta }(u) &=&\left[
\begin{array}{cc}
e_{\ \underline{i}}^{i\ }(u) & -N_{k}^{b}(u)e_{\ \underline{i}}^{k\ }(u) \\
0 & e_{\ \underline{a}}^{a\ }(u)%
\end{array}%
\right] ,  \label{vt2}
\end{eqnarray}%
we write equivalently the metric in the form
\begin{equation}
\mathbf{g}=\mathbf{g}_{\alpha \beta }\left( u\right) \mathbf{\vartheta }%
^{\alpha }\otimes \mathbf{\vartheta }^{\beta }=g_{ij}\left( u\right)
\vartheta ^{i}\otimes \vartheta ^{j}+h_{ab}\left( u\right) \ ^{\star
}\vartheta ^{a}\otimes \ ^{\star }\vartheta ^{b},  \label{metr}
\end{equation}%
where $g_{ij}\doteqdot \mathbf{g}\left( e_{i},e_{j}\right) $ and $%
h_{ab}\doteqdot \mathbf{g}\left( e_{a},e_{b}\right) $ \ and
\begin{equation*}
\mathbf{e}_{\alpha }=\mathbf{e}_{\alpha }^{\ \underline{\alpha }}\partial _{%
\underline{\alpha }}\mbox{ and }\mathbf{\vartheta }_{\ }^{\beta }=\mathbf{e}%
_{\ \underline{\beta }}^{\beta }du^{\underline{\beta }}.
\end{equation*}%
are vielbeins of type (\ref{dder}) and (\ref{ddif}) defined for arbitrary $%
N_{i}^{b}(u).$ We can consider a special class of manifolds provided with a
global splitting into conventional ''horizontal'' and ''vertical'' subspaces
(\ref{whitney}) induced by the ''off--diagonal'' terms $N_{i}^{b}(u)$ and
prescribed type of nonholonomic frame structure.

If the manifold $\mathbf{V}$ is (pseudo) Riemannian, there is a unique
linear connection (the Levi--Civita connection) $\nabla ,$ which is metric, $%
\nabla \mathbf{g=0,}$ and torsionless, $\ ^{\nabla }T=0.$ Nevertheless, the
connection $\nabla $ is not adapted to the nonintegrable distribution
induced by $N_{i}^{b}(u).$ In this case, 
for instance, in order to construct exact solutions parametrized by generic
off--diagonal metrics, or for investigating nonholonomic frame structures in
gravity models with nontrivial torsion, it is more convenient to work with
more general classes of linear connections which are N--adapted but contain
nontrivial torsion coefficients because of nontrivial nonholonomy
coefficients $W_{\alpha \beta }^{\gamma }$ (\ref{anhrel}).

For a splitting of a (pseudo) Riemannian--Cartan space of dimension $(n+m)$
(under certain constraints, we can consider (pseudo) Riemannian
configurations), the Lagrange and Finsler type geometries were modeled by
N--anholonomic structures as exact solutions of gravitational field
equations \cite{vex1,vt,vs1,dv}, see also Refs. \cite{v2,v3} for exact
solutions with nonmetricity.  One holds \cite{v1} the
\begin{result} \label{r6}
 The  geometry of any  Riemannian space of dimension $n+m$  where
 $n,m \geq 2$  (we can consider $n,m=1$ as special degenerated cases),
 provided with  off--diagonal metric structure of type (\ref{ansatz})
 can be equivalently modeled, by vielbein transforms of type
 (\ref{vt1}) and (\ref{vt2}) as a geometry of nonholonomic manifold
 enabled with N--connection structure $N_{i}^{b}(u)$ and 'more
 diagonalized' metric (\ref{metr}). 
\end{result}
For particular cases, we present the  
\begin{remark}
For certain special conditions
 when $n=m,$ $N_{i}^{b}=\ ^{(L)}N_{i}^{b}$ (\ref{cncl}) and the metric 
(\ref{metr}) is of type (\ref{slm}), a such Riemann space of even
 dimension is 'nonholonomically' equivalent to a Lagrange space (for
 the corresponding homogeneity conditions, see Definition \ref{dlargf},
 one obtains the equivalence to a Finsler space). 
\end{remark}
Roughly speaking, by prescribing corresponding nonholonomic frame
structures, we can model a Lagrange, or Finsler, geometry on 
a Riemannian manifold and, inversely, a Riemannian geometry is 'not
only a Riemannian one' but also could be a generalized Finsler
one. It is possible to define similar constructions for the
 (pseudo) Riemannian spaces. 
This is a quite surprising result if to compare it with the
"superficial" interpretation of the Finsler geometry as a nonlinear
extension, 'more sophisticate' on the tangent bundle, 
 of the Riemannian geometry. 

It is
known the fact that the first example of Finsler geometry was
considered in 1854 in the famous B. Riemann's hability thesis (see
historical  details and discussion 
in Refs. \cite{sh,bcs,ma2,v1}) who, for simplicity, restricted his
considerations only to the curvatures defined by quadratic
forms on hypersurfaces.
 Sure,  for B. Riemann, it  was unknown  the fact that if we
consider general (nonholonomic) frames with  associated nonlinear
connections (the E. Cartan's geometry, see Refs. in \cite{car1})
 and off--diagonal metrics, the
Finsler geometry may be derived naturally even from quadratic metric
forms being adapted to the N--connection structure.

 More rigorous geometric constructions involving the
Cartan--Miron metric connections and, respectively,
  the Berwald and Chern--Rund
 nonmetric  connections in Finsler geometry and generalizations, see
 more details in subsection \ref{sdlc},   
 result in equivalence theorems to certain types of 
 Riemann--Cartan nonholonomic
 manifolds (with nontirvial N--connection and torsion) and
 metric--affine nonholonomic manifolds (with additional nontrivial
 nonmetricity structures) \cite{v1}. 
  
This Result \ref{r6} give rise to an important:
\begin{conclusion} \label{ic} To study  generalized Finsler spinor and
  noncommutative geometries is necessary
  even if we restrict our considerations only to
  (non) commutative Riemannian geometries.
\end{conclusion}

For simplicity, in this work we
 restrict our considerations only to certain Riemannian commutative and
 noncommutative geometries when the N--connection and torsion are
 defined by corresponding nonholonomic frames.

\subsection{N--anholonomic manifolds} 
\label{snc}

Now we shall demonstrate how general N--connection structures define a
certain class of nonholonomic geometries. In this case, it is
convenient to work on a general manifold $\mathbf{V,}$ $\dim
\mathbf{V=}n+m,$ with global splitting, instead of the tangent bundle
$\widetilde{TM}.$ The constructions will contain  those from
geometric mechanics and gravity theories, as certain particular
cases.

Let $\mathbf{V}$ be a $(n+m)$--dimensional manifold. It is supposed that in
any point $u\in \mathbf{V}$ there is a local distribution (splitting)
 $\mathbf{V}_{u}=M_{u}\oplus V_{u},$
 where $M$ is a $n-$dimensional subspace and $V$ is
a $m$--dimensional subspace. The local coordinates (in
general, abstract ones both for holonomic and nonholonomic variables)
may be written in the
form $u=(x,y),$ or $u^{\alpha }=\left( x^{i},y^{a}\right).$  We denote by $\pi
^{\top }:T\mathbf{V}\rightarrow TM$ the differential of a map $\pi
:V^{n+m}\rightarrow V^{n}$ defined by fiber preserving morphisms of the
tangent bundles $T\mathbf{V}$ and $TM.$ The kernel of $\pi ^{\top }$ is just
the vertical subspace $v\mathbf{V}$ with a related inclusion mapping $i:v%
\mathbf{V}\rightarrow T\mathbf{V}.$

\begin{definition} \label{dncam}
A nonlinear connection (N--connection) $\mathbf{N}$ on a manifold $\mathbf{V}
$ is defined by the splitting on the left of an exact sequence
\begin{equation*}
0\rightarrow v\mathbf{V}\overset{i}{\rightarrow} T\mathbf{V}\rightarrow T%
\mathbf{V}/v\mathbf{V}\rightarrow 0,
\end{equation*}%
i. e. by a morphism of submanifolds $\mathbf{N:\ \ }T\mathbf{V}\rightarrow v%
\mathbf{V}$ such that $\mathbf{N\circ i}$ is the unity in $v\mathbf{V}.$
\end{definition}

In an equivalent form, we can say that a N--connection is defined by a
splitting to subspaces with a Whitney sum of conventional h--subspace, 
$\left( h\mathbf{V}\right) ,$ and v--subspace, $\left( v%
\mathbf{V}\right) ,$
\begin{equation}
T\mathbf{V}=h\mathbf{V}\oplus v\mathbf{V}  \label{whitney}
\end{equation}%
where $h\mathbf{V}$ is isomorphic to $M.$ This  generalizes  
the splitting considered in Result \ref{resws}.

Locally, a N--connection is defined by its coefficients $N_{i}^{a}(u),$%
\begin{equation}
\mathbf{N}=N_{i}^{a}(u)dx^{i}\otimes \frac{\partial }{\partial y^{a}}.
\label{nconcoef}
\end{equation}%
The well known class of linear connections consists a particular subclass
with the coefficients being linear on $y^{a},$ i. e. $N_{i}^{a}(u)=\Gamma
_{bj}^{a}(x)y^{b}.$

Any N--connection also defines a N--connection curvature
\begin{equation*}
\mathbf{\Omega }=\frac{1}{2}\Omega _{ij}^{a}d^{i}\wedge d^{j}\otimes
\partial _{a},
\end{equation*}%
with N--connection curvature coefficients given by formula (\ref{anhrel}). %

\begin{definition} \label{dnam}
A manifold \ $\mathbf{V}$ is called N--anholonomic if on the tangent space $T%
\mathbf{V}$ it is defined a local (nonintegrable) distribution (\ref{whitney}%
), i. e. $T\mathbf{V}$ is enabled with a N--connection (\ref{nconcoef})
inducing a vielbein structure (\ref{dder}) satisfying the nonholonomy
relations  (\ref{anhrel}) (such N--connections and
associated vielbeins may be general ones or any derived from a
Lagrange/ Finsler fundamental function).
\end{definition}

We note that the  boldfaced symbols are used for the spaces and
geometric objects provided/adapted to a N--connection structure. For
instance, a vector field $\mathbf{X}\in T\mathbf{V}$ \ is expressed $\mathbf{%
X}=(X \equiv\ ^{-}X,\ ^{\star }X),$ or $\mathbf{X}=X^{\alpha }\mathbf{e}_{\alpha
}=X^{i}e_{i}+X^{a}e_{a},$ where $X=\  ^{-}X = X^{i}e_{i}$ and $^{\star }X=X^{a}e_{a}$
state, respectively, the irreducible (adapted to the N--connection
structure) h-- and v--components of the vector (which
following Refs. \cite{ma1,ma2} is called a distinguished vectors, in brief,
d--vector). In a similar fashion, the geometric objects on $\mathbf{V}$ like
tensors, spinors, connections, ... are respectively defined and called 
 d--tensors, d--spinors, d--connections if they are adapted to the N--connection
splitting.\footnote{In order to emphasize h-- and v--splitting of any
d--objects $\mathbf{Y}, \mathbf{g},$ ... we shall write the
irreducible components as $\mathbf{Y} = (\ ^{-}Y,\ ^\star Y),$\ 
 $\mathbf{g} = (\ ^{-}g,\ ^\star g)$ but we shall omit "$-$" or
"$\star$" if the simplified denotations will not result in ambiguities.}  

\begin{definition} \label{ddms}
  A d--metric structure on N--anholonomic manifold
  $\mathbf{V}$ is defined by a symmetric d--tensor field of type
 $\mathbf{g}=[g,\ ^{\star}h]$ (\ref{metr}).
\end{definition}

 For any fixed values of
  coordinates $u=(x,y) \in \mathbf{V}$  a d--metric it defines a symmetric
 quadratic d--metric form,
\begin{equation} \label{qdmf}
 \mathbf{q}(\mathbf{x},\mathbf{y}) \doteq g_{ij} x^i x^j + h_{ab} y^a
 y^b,
\end{equation}
where the $n+m$--splitting is defined by the N--connection structure
and  $\mathbf{x} = x^i e_i + x^a e_a,\ \mathbf{y}=y^i e_i + y^a e_a
\in V^{n+m}.$

Any d--metric is parametrized by a generic off--diagonal matrix 
(\ref{ansatz}) if the coefficients are redefined with respect to a
 local coordinate basis (for corresponding parametrizations of the
 the data $[g,h,N]$ such ansatz model a  geometry of mechanics, or a
 Finsler like structure, in a  Riemann--Cartan--Weyl space provided
 with N--connection structure  \cite{v1,v2};
 for certain constraints, there are possible models derived as exact
solutions in Einsten gravity and noncommutative generalizations
 \cite{vex1,dv,v3}).

\begin{remark} There is a  special case when
 $\dim \mathbf{V=}n+n,\ h_{ab}\rightarrow g_{ij}$ and 
$N_{i}^{a}\rightarrow N_{\ i}^{j}$ in (\ref{metr}), which models
locally,  on $\mathbf{V,}$ 
 a tangent bundle  structure. We denote a
such space by $\widetilde{\mathbf{V}}_{(n,n)}.$ If the N--connection
and d--metric coefficients are just the canonical ones for the Lagrange
(Finsler) geometry (see, respectively, formulas (\ref{cncl}) and
(\ref{slm}) ),
 we model such locally anisotropic structures not
on a tangent bundle $TM$ but on a N--anholonomic manifold of dimension $2n.$ 
\end{remark}

We present some historical remarks on
 N--con\-nec\-ti\-ons and related subjects:\ The geometrical aspects of
the N--connection formalism has been studied since the first papers of E.
Cartan \cite{car1}\ and A. Kawaguchi \cite{ak1,ak2}\ (who used it in
component form for Finsler geometry). Then one should be mentioned the so
called Ehressman connection \cite{eh}) and the work of W. Barthel \cite{wb}
where the global definition of N--connection was given. In monographs \cite%
{ma1,ma2,mhss},  the N--connection formalism was elaborated in details
and applied to the geometry of generalized Finsler--Lagrange and
 Cartan--Hamilton spaces, see also the approaches \cite{ml1,ml2,fe}.It
 should be  noted that the works related to  nonholonomic geometry 
and N--connections have appeared many times in a rather dispersive way 
when different schools of authors from
geometry, mechanics and physics have worked many times not having relation
with another. We cite only some our recent results with explicit applications in
modern mathematical physics and particle and string theories: 
N--connection structures were modeled on Clifford and spinor bundles
\cite{vsp1,vsph,vsp2,vst},
 on superbundles and in some directions of (super) string theory %
\cite{vstr,vstr1}, as well in noncommutative geometry and gravity
 \cite{vnc,vncg,vncgs}.
 The idea to apply the N--connections formalism as a new geometric method
of constructing exact solutions in gravity theories was suggested in Refs. %
\cite{vex1,vkin} and developed in a number of works, see for instance, Ref. %
\cite{vt,vs1,dv}).

\section{Curvature of N--anholonomic Manifolds}

The geometry of nonholonomic manifolds has a long time history of yet
unfinished elaboration:\ For instance, in the review \cite{ver} it is
stated that it is probably impossible to construct an analog of the
Riemannian tensor for the general nonholonomic manifold.  In a more
recent review \cite{mironnh1}, it is emphasized that in the
past there were proposed well defined Riemannian tensors for a number
of spaces provided with nonholonomic distributions, like Finsler and
Lagrange spaces and various type of theirs higher order
generalizations, i. e. for nonholonomic manifolds possessing
corresponding  N--connection structures. As some examples of first
such  investigations, we cite the works  \cite{nhm1,nhm3,nhm4}. 
 In this section  we shall construct in explicit form the curvature tensor for 
 the N--anholonomic manifolds.

\subsection{Distinguished connections }
\label{sdlc}

On N--anholonomic manifolds, the geometric
 constructions can  be adapted to the N--connection structure:

\begin{definition} \label{ddc}
A distinguished connection (d--connection) $\mathbf{D}$ on a manifold $%
\mathbf{V}$ is a linear connection conserving under parallelism the Whitney
sum (\ref{whitney}) defining a general N--connection.
\end{definition}

The N--adapted components $\mathbf{\Gamma }_{\beta \gamma }^{\alpha }$ of a
d-connection $\mathbf{D}_{\alpha }=(\delta _{\alpha }\rfloor \mathbf{D})$
are defined by the equations $
\mathbf{D}_{\alpha }\delta _{\beta }=\mathbf{\Gamma }_{\ \alpha \beta
}^{\gamma }\delta _{\gamma },$
or
\begin{equation}
\mathbf{\Gamma }_{\ \alpha \beta }^{\gamma }\left( u\right) =\left( \mathbf{D%
}_{\alpha }\delta _{\beta }\right) \rfloor \delta ^{\gamma }.  \label{dcon1}
\end{equation}%
In its turn, this defines a N--adapted splitting into h-- and v--covariant
derivatives, $\mathbf{D}=D+\ ^{\star }D,$ where $D_{k}=\left(
L_{jk}^{i},L_{bk\;}^{a}\right) $ and $\ ^{\star }D_{c}=\left(
C_{jk}^{i},C_{bc}^{a}\right) $ are introduced as corresponding h- and
v--parametrizations of (\ref{dcon1}),%
\begin{equation*}
L_{jk}^{i}=\left( \mathbf{D}_{k}e_{j}\right) \rfloor \vartheta ^{i},\quad
L_{bk}^{a}=\left( \mathbf{D}_{k}e_{b}\right) \rfloor \vartheta
^{a},~C_{jc}^{i}=\left( \mathbf{D}_{c}e_{j}\right) \rfloor \vartheta
^{i},\quad C_{bc}^{a}=\left( \mathbf{D}_{c}e_{b}\right) \rfloor \vartheta
^{a}.
\end{equation*}%
The components $\mathbf{\Gamma }_{\ \alpha \beta }^{\gamma }=\left(
L_{jk}^{i},L_{bk}^{a},C_{jc}^{i},C_{bc}^{a}\right) $ completely define a
d--connection $\mathbf{D}$ on a N--anholonomic manifold $\mathbf{V}.$

The simplest way to perform computations with d--connections is to use
N--adapted differential forms like $\mathbf{\Gamma }_{\beta }^{\alpha }=%
\mathbf{\Gamma }_{\beta \gamma }^{\alpha }\mathbf{\vartheta }^{\gamma }$
with the coefficients defined with respect to N--elongate bases (\ref{ddif}) and (\ref{dder}).

The torsion  of  d--connection $\mathbf{D}$ is defined by the
usual formula%
\begin{equation*}
\mathbf{T}(\mathbf{X},\mathbf{Y})\doteqdot \mathbf{D}_{X}\mathbf{D}%
_{Y}-\mathbf{D}_{Y}\mathbf{D}_{X} -[\mathbf{X},\mathbf{Y}].
\end{equation*}%

\begin{theorem}
The torsion $\mathbf{T}^{\alpha }\doteqdot \mathbf{D\vartheta }^{\alpha }=d%
\mathbf{\vartheta }^{\alpha }+\Gamma _{\beta }^{\alpha }\wedge \mathbf{%
\vartheta }^{\beta }$ of a d--connection has the irreducible h- v--
components (d--torsions)with N--adapted coefficients
\begin{eqnarray}
T_{\ jk}^{i} &=&L_{\ [jk]}^{i},\ T_{\ ja}^{i}=-T_{\ aj}^{i}=C_{\ ja}^{i},\
T_{\ ji}^{a}=\Omega _{\ ji}^{a},\   \notag \\
T_{\ bi}^{a} &=&T_{\ ib}^{a}=\frac{\partial N_{i}^{a}}{\partial y^{b}}-L_{\
bi}^{a},\ T_{\ bc}^{a}=C_{\ [bc]}^{a}.  \label{dtors}
\end{eqnarray}
\end{theorem}

\begin{proof}
By a straightforward calculation we can verify the formulas.
\end{proof}

The Levi--Civita linear connection $\nabla =\{^{\nabla }\mathbf{\Gamma }%
_{\beta \gamma }^{\alpha }\},$ with vanishing both torsion and
nonmetricity, is not adapted to the global splitting
(\ref{whitney}). One holds:

\begin{proposition} \label{pcnc}
There is a preferred, canonical d--connection structure, 
$\widehat{\mathbf{D}},$ on N--anholonomic manifold $\mathbf{V}$ constructed only
from the metric and N--con\-nec\-ti\-on coefficients $%
[g_{ij},h_{ab},N_{i}^{a}]$ and satisfying the metricity conditions $\widehat{\mathbf{D}%
}\mathbf{g}=0$ and $\widehat{T}_{\ jk}^{i}=0$ and $\widehat{T}_{\ bc}^{a}=0.$
\end{proposition}

\begin{proof}
By straightforward calculations with respect to the N--adapted bases (\ref%
{ddif}) and (\ref{dder}), we can verify that the connection
\begin{equation}
\widehat{\mathbf{\Gamma }}_{\beta \gamma }^{\alpha }=\ ^{\nabla }\mathbf{%
\Gamma }_{\beta \gamma }^{\alpha }+\ \widehat{\mathbf{P}}_{\beta \gamma
}^{\alpha }  \label{cdc}
\end{equation}%
with the deformation d--tensor
\begin{equation*}
\widehat{\mathbf{P}}_{\beta \gamma }^{\alpha }=(P_{jk}^{i}=0,P_{bk}^{a}=%
\frac{\partial N_{k}^{a}}{\partial y^{b}},P_{jc}^{i}=-\frac{1}{2}%
g^{ik}\Omega _{\ kj}^{a}h_{ca},P_{bc}^{a}=0)
\end{equation*}%
satisfies the conditions of this Proposition. It should be noted that, in
general, the components $\widehat{T}_{\ ja}^{i},\ \widehat{T}_{\ ji}^{a}$
and $\widehat{T}_{\ bi}^{a}$ are not zero. This is an anholonomic frame (or,
equivalently, off--diagonal metric) effect.
\end{proof}

With respect to the N--adapted frames, the coefficients\newline
$\widehat{\mathbf{\Gamma }}_{\ \alpha \beta }^{\gamma }=\left( \widehat{L}%
_{jk}^{i},\widehat{L}_{bk}^{a},\widehat{C}_{jc}^{i},\widehat{C}%
_{bc}^{a}\right) $ are computed:
\begin{eqnarray}
\widehat{L}_{jk}^{i} &=&\frac{1}{2}g^{ir}\left( \frac{\delta g_{jr}}{%
\partial x^{k}}+\frac{\delta g_{kr}}{\partial x^{j}}-\frac{\delta g_{jk}}{%
\partial x^{r}}\right) ,  \label{candcon} \\
\widehat{L}_{bk}^{a} &=&\frac{\partial N_{k}^{a}}{\partial y^{b}}+\frac{1}{2}%
h^{ac}\left( \frac{\delta h_{bc}}{\partial x^{k}}-\frac{\partial N_{k}^{d}}{%
\partial y^{b}}h_{dc}-\frac{\partial N_{k}^{d}}{\partial y^{c}}h_{db}\right)
,  \notag \\
\widehat{C}_{jc}^{i} &=&\frac{1}{2}g^{ik}\frac{\partial g_{jk}}{\partial
y^{c}},  \notag \\
\widehat{C}_{bc}^{a} &=&\frac{1}{2}h^{ad}\left( \frac{\partial h_{bd}}{%
\partial y^{c}}+\frac{\partial h_{cd}}{\partial y^{b}}-\frac{\partial h_{bc}%
}{\partial y^{d}}\right) .  \notag
\end{eqnarray}%
The d--connection (\ref{candcon}) defines the 'most minimal'    N--adapted
 extension of the Levi--Civita connection in order to preserve the
 metricity condition and to have zero torsions on the h-- and
 v--subspaces (the rest of nonzero torsion coefficients are defined by
 the condition of compatibility with the N--connection splitting). 

\begin{remark}
The canonical d--connection $\widehat{\mathbf{D}}$ (\ref{candcon}) for a
local modelling of a $\widetilde{TM}$ space on $\widetilde{\mathbf{V}}%
_{(n,n)}$ is defined by the coefficients
 $\widehat{\mathbf{\Gamma }}_{\ \alpha \beta
}^{\gamma }=(\widehat{L}_{jk}^{i},\widehat{C}_{jk}^{i})\ $ with%
\begin{equation}
\widehat{L}_{jk}^{i}=\frac{1}{2}g^{ir}\left( \frac{\delta g_{jr}}{\partial
x^{k}}+\frac{\delta g_{kr}}{\partial x^{j}}-\frac{\delta g_{jk}}{\partial
x^{r}}\right) , \widehat{C}_{jk}^{i}=\frac{1}{2}g^{ir}\left( \frac{\partial
g_{jr}}{\partial y^{k}}+\frac{\partial g_{kr}}{\partial y^{j}}-\frac{%
\partial g_{jk}}{\partial y^{r}}\right)  \label{candcon1} 
\end{equation}%
 computed with respect to N--adapted bases (\ref{dder}) and
 (\ref{ddif}) when  $\widehat{L}_{jk}^{i}$ and $\widehat{C}_{jk}^{i}$ define
respectively the canonical h-- and v--covariant derivations.
\end{remark}

Various models of Finsler geometry and generalizations were elaborated
by using different types of d--connections which satisfy, or not, the
compatibility conditions with a fixed d--metric structure (for
instance, with  a Sasaki type one). Let us consider the main examples:
\begin{example} The Cartan's d--connection \cite{car1}  with 
  the coefficients (\ref{candcon1})  was defined by 
   some generalized Christoffel symbols with the aim to have a 'minimal'
  torsion and to preserve the metricity condition. This approach was
  developed for generalized Lagrange spaces and on vector bundles
  provided with N--connection structure \cite{ma1,ma2}  by introducing the
  canonical d--connection (\ref{candcon}). The direction emphasized
  metric compatible and N--adapted geometric constructions. 
\end{example}

An alternative class of Finsler geometries is concluded in monographs
\cite{bcs,sh}: 
\begin{example} \label{ecrdc} It was the idea of C. C. Chern
 \cite{chern} (latter also
proposed by H. Rund \cite{rund}) to consider a d--connection
 $\ ^{[Chern]}{\mathbf{\Gamma }}_{\ \alpha \beta
}^{\gamma }=(\widehat{L}_{jk}^{i},{C}_{jk}^{i}= 0)\ $ and to work not
on a tangent bundle $TM$ but to try to 'keep maximally' the constructions
 on the base manifold $M.$ The Chern d--connection, as well
the Berwald d--connection   $\ ^{[Berwald]}{\mathbf{\Gamma }}_{\ \alpha \beta
}^{\gamma }=({L}_{jk}^{i}=\frac{\partial {N}_{k}^{i}}{\partial y^j}
,{C}_{jk}^{i}= 0)\ $ \cite{bw}, are not subjected to the metricity
conditions. 
\end{example} 
We note that the constructions mentioned in the last example 
 define certain  'nonmetric geometries' 
(a Finsler modification of the Riemann--Cartan--Weyl
spaces). For the Chern's connection, the torsion vanishes but there is
a nontrivial nonmetricity. A detailed study and classification  of Finsler--affine spaces
with general  nontrivial N--connection, torsion and nonmetricity was
recently performed in Refs. \cite{v1,v2,v3}. 
Here we also note that we may consider any  linear connection can be
generated by deformations of type 
\begin{equation}
\mathbf{\Gamma }_{\beta \gamma }^{\alpha
}=\widehat{\mathbf{\Gamma }}_{\beta \gamma }^{\alpha }+\mathbf{P}_{\beta
\gamma }^{\alpha }. \label{defcon}
\end{equation}
This  splits all geometric objects into canonical and
post-canonical pieces which results in N--adapted geometric
constructions.

In order to define spinors on generalized Lagrange and
Finsler spaces \cite{vsp1,vsph,vst,vsp2} the canonical d--connection
and Cartan's d--connection were used because the metric compatibility
allows the simplest definition of Clifford structures locally adapted
to the N--connection. This is also the simplest way to define
the Dirac operator for generalized Finsler spaces and to extend the
constructions to noncommutative Finsler geometry
\cite{vnc,vncg,vncgs}. The geometric constructions with general metric
compatible affine connection (with torsion) are preferred in
modern gravity and string theories. Nevertheless, the geometrical and
physical models with generic nonmetricity also present certain
interest \cite{hehl,v1,v2,v3} (see also \cite{maj} where nonmetricity
is considered to be important in quantum group co gravity). In such
cases, we can use deformations of connection (\ref{defcon}) in order
to 'deform', for instance, the spinorial geometric constructions
defined by the canonical d--connection and to transform them into 
certain 'nonmetric' configurations.

\subsection{ Curvature of  d--connections}

The curvature of a  d--connection $\mathbf{D}$ on an N--anholonomic
manifold is defined by the usual formula%
\begin{equation*}
\mathbf{R}(\mathbf{X},\mathbf{Y})\mathbf{Z}\doteqdot \mathbf{D}_{X}\mathbf{D}%
_{Y}\mathbf{Z}-\mathbf{D}_{Y}\mathbf{D}_{X}\mathbf{Z-D}_{[X,X]}\mathbf{Z.}
\end{equation*}%

By straightforward calculations we prove:
\begin{theorem}
The curvature $\mathcal{R}_{\ \beta }^{\alpha }\doteqdot \mathbf{D\Gamma }%
_{\beta }^{\alpha }=d\mathbf{\Gamma }_{\beta }^{\alpha }-\mathbf{\Gamma }%
_{\beta }^{\gamma }\wedge \mathbf{\Gamma }_{\gamma }^{\alpha }$ of a
d--connection $\mathcal{D} \doteq \mathbf{\Gamma }_{\gamma }^{\alpha }$ has the irreducible h-
v-- components (d--curvatures) of $\mathbf{R}_{\ \beta \gamma \delta
}^{\alpha }$,%
\begin{eqnarray}
R_{\ hjk}^{i} &=&e_{k}L_{\ hj}^{i}-e_{j}L_{\ hk}^{i}+L_{\ hj}^{m}L_{\
mk}^{i}-L_{\ hk}^{m}L_{\ mj}^{i}-C_{\ ha}^{i}\Omega _{\ kj}^{a},  \notag \\
R_{\ bjk}^{a} &=&e_{k}L_{\ bj}^{a}-e_{j}L_{\ bk}^{a}+L_{\ bj}^{c}L_{\
ck}^{a}-L_{\ bk}^{c}L_{\ cj}^{a}-C_{\ bc}^{a}\Omega _{\ kj}^{c},  \notag \\
R_{\ jka}^{i} &=&e_{a}L_{\ jk}^{i}-D_{k}C_{\ ja}^{i}+C_{\ jb}^{i}T_{\
ka}^{b},  \label{dcurv} \\
R_{\ bka}^{c} &=&e_{a}L_{\ bk}^{c}-D_{k}C_{\ ba}^{c}+C_{\ bd}^{c}T_{\
ka}^{c},  \notag \\
R_{\ jbc}^{i} &=&e_{c}C_{\ jb}^{i}-e_{b}C_{\ jc}^{i}+C_{\ jb}^{h}C_{\
hc}^{i}-C_{\ jc}^{h}C_{\ hb}^{i},  \notag \\
R_{\ bcd}^{a} &=&e_{d}C_{\ bc}^{a}-e_{c}C_{\ bd}^{a}+C_{\ bc}^{e}C_{\
ed}^{a}-C_{\ bd}^{e}C_{\ ec}^{a}.  \notag
\end{eqnarray}
\end{theorem}

\begin{remark}
For an N--anholonomic manifold $\widetilde{\mathbf{V}}_{(n,n)}$ provided
with N--sym\-plet\-ic canonical d--connection $\widehat{\mathbf{\Gamma }}_{\
\alpha \beta }^{\tau },$ see (\ref{candcon1}), the d--curvatures (\ref{dcurv}%
) reduces to three irreducible components
\begin{eqnarray}
R_{\ hjk}^{i} &=&e_{k}L_{\ hj}^{i}-e_{j}L_{\ hk}^{i}+L_{\ hj}^{m}L_{\
mk}^{i}-L_{\ hk}^{m}L_{\ mj}^{i}-C_{\ ha}^{i}\Omega _{\ kj}^{a},  \notag \\
R_{\ jka}^{i} &=&e_{a}L_{\ jk}^{i}-D_{k}C_{\ ja}^{i}+C_{\ jb}^{i}T_{\
ka}^{b},  \label{dcurv1} \\
R_{\ bcd}^{a} &=&e_{d}C_{\ bc}^{a}-e_{c}C_{\ bd}^{a}+C_{\ bc}^{e}C_{\
ed}^{a}-C_{\ bd}^{e}C_{\ ec}^{a}  \notag
\end{eqnarray}%
where all indices $i,j,k...$ and $a,b,..$ run the same values but label the
components with respect to different h-- or v--frames.
\end{remark}

Contracting respectively the components of (\ref{dcurv}) and (\ref{dcurv1})
we prove:

\begin{corollary}
The Ricci d--tensor $\mathbf{R}_{\alpha \beta }\doteqdot \mathbf{R}_{\
\alpha \beta \tau }^{\tau }$ has the irreducible h- v--components%
\begin{equation}
R_{ij}\doteqdot R_{\ ijk}^{k},\ \ R_{ia}\doteqdot -R_{\ ika}^{k},\
R_{ai}\doteqdot R_{\ aib}^{b},\ R_{ab}\doteqdot R_{\ abc}^{c},
\label{dricci}
\end{equation}%
for a general N--holonomic manifold $\mathbf{V,}$ and
\begin{equation}
R_{ij}\doteqdot R_{\ ijk}^{k},\ \ R_{ia}\doteqdot -R_{\ ika}^{k},\ \
R_{ab}\doteqdot R_{\ abc}^{c},  \label{dricci1}
\end{equation}%
for an N--anholonomic manifold $\widetilde{\mathbf{V}}_{(n,n)}.$
\end{corollary}

\begin{corollary}
The scalar curvature of a d--connection is
\begin{equation} \label{sdccurv}
\overleftarrow{\mathbf{R}} \doteqdot \mathbf{g}^{\alpha \beta }\mathbf{R}%
_{\alpha \beta }=g^{ij}R_{ij}+h^{ab}R_{ab}, 
\end{equation}
defined by the "pure" h-- and v--components of (\ref{dricci1}).
\end{corollary}

\begin{corollary}
The Einstein d--densor is computed $\mathbf{G}_{\alpha \beta }=\mathbf{R}%
_{\alpha \beta }-\frac{1}{2}\mathbf{g}_{\alpha \beta }\overleftarrow{\mathbf{%
R}}.$
\end{corollary}

For physical applications, the Riemann, Ricci and Einstein d--tensors
can be computed for the canonical d--connection. We can 
redefine the constructions for arbitrary d--connections by using the
corresponding deformation tensors like in (\ref{defcon}), for instance, 
\begin{equation}
\mathcal{R}_{\ \beta }^{\alpha }=\widehat{\mathcal{R}}_{\ \beta }^{\alpha }+%
\mathbf{D}\mathcal{P}_{\ \beta }^{\alpha }+\mathcal{P}_{\ \gamma }^{\alpha
}\wedge \mathcal{P}_{\ \beta }^{\gamma }  \label{deformcurv}
\end{equation}%
for $\mathcal{P}_{\beta }^{\alpha }=\mathbf{P}_{\beta \gamma }^{\alpha }%
\mathbf{\vartheta }^{\gamma }.$ 
A set of examples of such deformations are analyzed in
Refs. \cite{v1,v2,v3}.

\section{Noncommutative N--Anholonomic Spaces}

In this section, we outline how the analogs of basic objects in 
commutative geometry of N--anholonomic manifolds, such as
vector/tangent bundles, N-- and d--connections can be defined in
noncommutative geometry \cite{vncg,vncgs}. We note that the
 A. Connes' functional analytic approach \cite{connes1} to the
noncommutative topology and geometry is based on the theory of
noncommutative $C^{\ast }$--algebras. Any commutative $C^{\ast }$--algebra
can be realized as the $C^{\ast }$--algebra of complex valued functions over
locally compact Hausdorff space. A noncommutative $C^{\ast }$--algebra can
be thought of as the algebra of continuous functions on some 'noncommutative
space' (see main definitions and results in Refs.  \cite{connes1,bondia,landi,madore}).

The starting idea of noncommutative geometry is to derive the
geometric properties of ``commutative'' spaces from their algebras of
functions  characterized by involutive algebras of operators by
dropping the condition of commutativity (see the Gelfand and Naimark
theorem \cite{gelfand}). A space topology is defined by the algebra of
commutative continuous  functions, but the geometric constructions request a
differentiable structure. Usually, one considers a differentiable and
compact manifold $M,\ dim N=n$ (there is an open problem how to
include in noncommutative geometry spaces with indefinite metric signature like
pseudo--Euclidean and pseudo--Riemannian ones). In order to construct
models of commutative and noncommutative differential geometries it is
more or less obvious that the class of algebras of smooth functions,
 $\mathcal{C}\doteq C^{\infty}(M)$ is more appropriate. If $M$ is a
 smooth manifold, it is possible to reconstruct this manifold with its
 smooth structure and the attached objects (differential forms,
 etc...) by starting from  $\mathcal{C}$ considered as an abstract
 (commutative) unity $*$--algebra with involution. As a set $M$ can
 be identified with the set of characters of $\mathcal{C},$ but its
 differential structure is connected with the abundance of derivations
 of $\mathcal{C}$ which identify with the smooth vector fields on $M.$
There are two standard constructions: 1) when the vector fields are
 considered to be the derivations of  $\mathcal{C}$ (into itself) or
 2) one considers a generalization of the calculus of differential
 forms which is the Kahler differential calculus (see, details in
 Lectures \cite{dubois}). The noncommutative versions of differential
 geometry may be elaborated if the algebra of smooth complex functions
 on a smooth manifold is replaced by a noncommutative associative
 unity complex $*$--algebra $\mathcal{A}.$

The geometry of commutative gauge and gravity theories is derived
from the notions of connections (linear and nonlinear ones),
 metrics and frames of references on manifolds and
vector bundle spaces. The possibility of extending such theories to some
noncommutative models is based on the Serre--Swan theorem \cite{swan}
stating that there is a complete equivalence between the category of
(smooth) vector bundles over a smooth compact space (with bundle maps) and
the category of porjective modules of finite type over commutative algebras
and module morphisms. So, the space $\Gamma \left(
E\right) $ of smooth sections of a vector bundle $E$ over a compact space is
a projective module of finite type over the algebra $C\left( M\right) $ of
smooth functions over $M$ and any finite projective $C\left( M\right) $%
--module can be realized as the module of sections of some vector bundle
over $M.$ This construction may be extended if a noncommutative algebra $%
\mathcal{A}$ is taken as the starting ingredient:\ the noncommutative
analogue of vector bundles are projective modules of finite type over $%
\mathcal{A}$. This way one developed a theory of linear connections which
culminates in the definition of Yang--Mills type actions or, by some much
more general settings, one reproduced Lagrangians for the Standard model
with its Higgs sector or different type of gravity and Kaluza--Klein models
(see, for instance, Refs \cite{connes1,madore}).

\subsection{Modules as bundles}

A vector space $\mathcal{E}$ over the complex number field $\C$ can be
defined also as a right module of an algebra $\mathcal{A}$ over $\C$ \ which
carries a right representation of $\mathcal{A},$ when for every map of
elements $\mathcal{E}$ $\times \mathcal{A}\ni \left( \eta ,a\right)
\rightarrow \eta a\in \mathcal{E}$ one hold the properties 
\begin{equation*}
\lambda (ab)=(\lambda a)b,~\lambda (a+b)=\lambda a+\lambda b,~(\lambda +\mu
)a=\lambda a+\mu a
\end{equation*}%
fro every $\lambda ,\mu \in \mathcal{E}$ and $a,b\in \mathcal{A}.$

Having two $\mathcal{A}$--modules $\mathcal{E}$ and $\mathcal{F},$ a
morphism of $\mathcal{E}$ into $\mathcal{F}$ is \ any linear map $\rho :%
\mathcal{E}$ $\rightarrow $ $\mathcal{F}$  which is also $\mathcal{A}$%
--linear, i. e. $\rho (\eta a)=\rho (\eta )a$ for every $\eta \in \mathcal{E}
$ and $a\in \mathcal{A}.$

We can define in a similar (dual) manner the left modules and theirs
morphisms which are distinct from the right ones for noncommutative algebras 
$\mathcal{A}.$ A bimodule over an algebra $\mathcal{A}$ is a vector space $%
\mathcal{E}$ which carries both a left and right module
structures. The bimodule structure is important for modeling of real
geometries starting from complex structures. We may
define the opposite algebra $\mathcal{A}^{o}$ with elements $a^{o}$ being in
bijective correspondence with the elements \ $a\in \mathcal{A}$ while the
multiplication is given by $\mathcal{\,}a^{o}b^{o}=\left( ba\right) ^{o}.$A
right (respectively, left) $\mathcal{A}$--module $\mathcal{E}$ is connected
to a left (respectively right) $\mathcal{A}^{o}$--module via relations $%
a^{o}\eta =\eta a^{o}$ (respectively, $a\eta =\eta a).$ 
One introduces the enveloping algebra $\mathcal{A}^{\varepsilon }=\mathcal{A}%
\otimes _{\C}\mathcal{A}^{o};$ any $\mathcal{A}$--bimodule $\mathcal{E}$ can
be regarded as a right [left] $\mathcal{A}^{\varepsilon }$--module by
setting $\eta \left( a\otimes b^{o}\right) =b\eta a$ $\quad \left[ \left(
a\otimes b^{o}\right) \eta =a\eta b\right] .$

For a (for instance, right) module $\mathcal{E}$ , we may introduce a family
of elements $\left( e_{t}\right) _{t\in T}$ parametrized by any (finite or
infinite) directed set $T$ for which any element $\eta \in \mathcal{E}$ is
expressed as a combination (in general, in \ more than one manner) $\eta
=\sum\nolimits_{t\in T}e_{t}a_{t}$ with $a_{t}\in \mathcal{A}$ and only a
finite number of non vanishing terms in the sum. A family $\left(
e_{t}\right) _{t\in T}$ is free if it consists from linearly independent
elements and defines a basis if any element $\eta \in \mathcal{E}$ can be
written as a unique combination (sum). One says a module to be free if it
admits a basis. The module $\mathcal{E}$ is said to be of finite type if \
it is finitely generated, i. e. it admits a generating family of finite
cardinality.

Let us consider the module $\mathcal{A}^{K}\doteqdot \C^{K}\otimes _{\C}%
\mathcal{A}.$ The elements of this module can be thought as $K$--dimensional
vectors with entries in $\mathcal{A}$ and written uniquely as a linear
combination $\eta =\sum\nolimits_{t=1}^{K}e_{t}a_{t}$ were the basis $e_{t}$
identified with the canonical basis of $\C^{K}.$ This is a free and finite
type module. In general, we can have bases of different cardinality.
However, if a module $\mathcal{E}$ \ is of finite type there is always an
integer $K$ and a module surjection $\rho :\mathcal{A}^{K}\rightarrow 
\mathcal{E}$ with a base being a image of a free basis, $\epsilon _{j}=\rho
(e_{j});j=1,2,...,K.$

We say that a right $\mathcal{A}$--module $\mathcal{E}$ is projective if for
every surjective module morphism $\rho :\mathcal{M}$ $\rightarrow $ $%
\mathcal{N}$ splits, i. e. there exists a module morphism \ $s:\mathcal{E}$ $%
\rightarrow $ $\mathcal{M}$ such that $\rho \circ s=id_{\mathcal{E}}.$ There
are different definitions of porjective modules (see Ref. \cite{landi} on
properties of such modules). Here we note the property that if a $\mathcal{A}
$--module $\mathcal{E}$ is projective, there exists a free module $\mathcal{F%
}$ and a module $\mathcal{E}^{\prime }$ (being a priory projective) such
that $\mathcal{F}=\mathcal{E}\oplus \mathcal{E}^{\prime }.$

For the right $\mathcal{A}$--module $\mathcal{E}$ being projective and of
finite type with surjection $\rho :\mathcal{A}^{K}\rightarrow \mathcal{E}$
and following the projective property we can find a lift $\widetilde{\lambda 
}:\mathcal{E}$ $\rightarrow $ $\mathcal{A}^{K}$ such that $\rho \circ 
\widetilde{\lambda }=id_{\mathcal{E}}.$ There is a proof of the property
that the module $\mathcal{E}$ is projective of finite type over $\mathcal{A}$
if and only if there exists an idempotent $p\in End_{\mathcal{A}}\mathcal{A}%
^{K}=M_{K}(\mathcal{A}),$ $p^{2}=p,$ the $M_{K}(\mathcal{A})$ denoting the
algebra of $K\times K$ matrices with entry in $\mathcal{A},$ such that $%
\mathcal{E}=p\mathcal{A}^{K}.$ We may associate the elements of $\mathcal{E}$
to $K$--dimensional column vectors whose elements are in $\mathcal{A},$ the
collection \ of which are invariant under the map $p,$ $\mathcal{E}$ $=\{\xi
=(\xi _{1},...,\xi _{K});\xi _{j}\in \mathcal{A},~p\xi =\xi \}.$ For
simplicity, we shall use the term finite projective to mean projective of
finite type.

\subsection{Nonlinear connections in projective modules}

The nonlinear connection (N--connection) for noncommutative spaces
can be defined similarly to commutative spaces by considering instead of
usual vector bundles theirs noncommutative analogs defined as finite
projective modules over noncommutative algebras \cite{vncg}.
 The explicit constructions depend on the type of differential
 calculus we use for definition of tangent structures and theirs
 maps. In this subsection, we shall consider such projective modules
  provided with N--connection which define noncommutative analogous
  both of vector bundles and of N--anholonomic manifolds (see
  Definition \ref{dnam}).

In general, one can be defined several differential calculi over a
given algebra $\mathcal{A}$ 
 (for a more detailed discussion within the context of
noncommutative geometry, see Refs. \cite{connes1,madore}).\ For
simplicity, in this work we consider that a differential calculus on
$\mathcal{A}$ is fixed,  which means that we choose a (graded) algebra
 $\Omega ^{\ast }(\mathcal{A})=\cup _{p}\Omega ^{p}(\mathcal{A})$
  giving a differential structure to $\mathcal{A}.$ The elements
 of  $\Omega ^{p}(\mathcal{A})$ are called $p$--forms. There is a linear map $d$
which takes $p$--forms into $(p+1)$--forms and which satisfies a graded
Leibniz rule as well the condition $d^{2}=0.$ 
By definition $\Omega ^{0}(\mathcal{A})=\mathcal{A}.$

The differential $df$ of a real or complex variable on a
N--anholonomic manifold  $\mathbf{V}$  
\begin{eqnarray*}
df &=&\delta _{i}f~dx^{i}+\partial _{a}f~\delta y^{a}, \\
\delta _{i}f~ &=&\partial _{i}f-N_{i}^{a}\partial _{a}f~,~\delta
y^{a}=dy^{a}+N_{i}^{a}dx^{i}, 
\end{eqnarray*}%
where the N--elongated derivatives and differentials are defined
respectively  by formulas (\ref{dder}) and (\ref{ddif}),  
in the noncommutative case is replaced by a distinguished commutator
(d--commutator)%
\begin{equation*}
\overline{d}f=\left[ F,f\right] =\left[ F^{[h]},f\right] +\left[ F^{[v]},f%
\right]
\end{equation*}%
where the operator $F^{[h]}$ $\ (F^{[v]})$ acts on the horizontal
(vertical) projective submodule and this operator is  defined by a
fixed differential calculus $\Omega ^{\ast }(\mathcal{A}^{[h]})$
($\Omega ^{\ast }(\mathcal{A}^{[v]}))$ on the so--called horizontal
(vertical)  $\mathcal{A}^{[h]}$ ($\mathcal{A}^{[v]})$ algebras. We
conclude that in order to elaborated noncommutative versions of
N--anholonomic manifolds we need couples of 'horizontal' and 'vertical'
operators which reflects the nonholonomic splitting given by the
N--connection structure. 

Let us consider instead of a N--anholonomic manifold $\mathbf{V}$ 
 an $\mathcal{A}$--module $\mathcal{E}$ being projective and of finite
 type.  For a fixed differential calculus on $\mathcal{E}$ we define
 the  tangent structures $T\mathcal{E}.$  
\begin{definition} \label{dncamnc}
A nonlinear connection (N--connection) $\mathbf{N}$ on  an $\mathcal{A}$%
--module $\mathcal{E}$  is defined by the splitting on the left of an
 exact sequence of finite projective $\mathcal{A}$--moduli 
\begin{equation*}
0\rightarrow v\mathcal{E}\overset{i}{\rightarrow} T\mathcal{E}\rightarrow T%
\mathcal{E}/v\mathcal{E}\rightarrow 0,
\end{equation*}%
i. e. by a morphism of submanifolds $\mathbf{N:\ }T\mathcal{E}\rightarrow v%
\mathcal{E}$ such that $\mathbf{N\circ i}$ is the unity in $v\mathcal{E}.$
\end{definition}
In an equivalent form, we can say that a N--connection is defined by a
splitting to projective submodules with a Whitney sum of conventional 
 h--submodule, $\left( h\mathcal{E}\right) ,$ and v--submodule,
 $\left( v\mathcal{E}\right) ,$
\begin{equation}
T\mathcal{E}=h\mathcal{E}\oplus v\mathcal{E}.  \label{whitneync}
\end{equation}%
We note that locally $h\mathcal{E}$ is isomorphic to $TM$ where $M$ is a
 differential compact manifold of dimension $n.$

The Definition \ref{dncamnc} reconsiders  for noncommutative spaces the
Definition \ref{dncam}. In result, we may generalize the concept of
'commutative' N--anholonomic space:

\begin{definition}  A N--anholonomic noncommutative space $\mathcal{E}_N$ is an 
 $\mathcal{A}$--mo\-du\-le $\mathcal{E}$  possessing a  tangent structure
 $T\mathcal{E}$   defined  by a Whitney sum of projective submodules
 (\ref{whitneync}).
\end{definition}
Such geometric  constructions  depend  on the type of fixed differential
calculus, i. e. on the procedure how the tangent spaces are
defined.

\begin{remark}
Locally always N--connections exist, but it is not obvious if they
could be glued together. In the classical case of vector bundles over
paracompact manifolds this is possible \cite{ma1}. If there is an
appropriate partition of unity, a similar result can be proved for
finite projective modules. 
 For certain applications,  it is more convenient to use the
Dirac operator already defined on N--anholonomic manifolds, see
Section \ref{sdonst}. 
\end{remark}

In order to understand how the N--connection structure may be taken
 into account on noncommutative spaces but distinguished from the
 class of  linear  gauge fields,  we analyze   an example:

\subsection{Commutative and noncommutative gauge d--fields}

Let us consider a N--anholonomic manifold $\mathbf{V}$  and a 
 vector bundle $\beta =\left( B,\pi ,\mathbf{V}\right) $ with
 $\pi :B\rightarrow \mathbf{V}$
with a typical $k$-dimensional vector fiber. In local coordinates a linear
connection (a gauge field) in $\beta $ is given by a collection of
differential operators%
\begin{equation*}
\nabla _{\alpha }=D_{\alpha }+B_{\alpha }(u),
\end{equation*}%
acting on $T\xi _{N}$ where 
\begin{equation*}
D_{\alpha }=\delta _{\alpha }\pm \Gamma _{\cdot \alpha }^{\cdot }\mbox{ with
}D_{i}=\delta _{i}\pm \Gamma _{\cdot i}^{\cdot }\mbox{ and }D_{a}=\partial
_{a}\pm \Gamma _{\cdot a}^{\cdot }
\end{equation*}
is a d--connection in $\mathbf{V}$  ($\alpha =1,2,...,n+m),$ with the
operator  $\delta _{\alpha },$
 being  N--elongated as in (\ref{dder}), $u=(x,y)\in \xi _{N}$ and $%
B_{\alpha }$ are $k\times k$--matrix valued functions. For every vector
field 
\begin{equation*}
X=X^{\alpha }(u)\delta _{\alpha }=X^{i}(u)\delta _{i}+X^{a}(u)\partial
_{a}\in T\mathbf{V} 
\end{equation*}
we can consider the operator 
\begin{equation}
X^{\alpha }(u)\nabla _{\alpha }(f\cdot s)=f\cdot \nabla
_{X}s+\delta _{X}f\cdot s  \label{rul1c}
\end{equation}%
for any section $s\in \mathcal{B}$ \ and function
 $f\in C^{\infty }(\mathbf{V}),$ where%
\begin{equation*}
\delta _{X}f=X^{\alpha }\delta _{\alpha }~\mbox{ and
}\nabla _{fX}=f\nabla _{X}.
\end{equation*}%
In the simplest definition we assume that there is a Lie algebra $\mathcal{GL%
}B$ that acts on associative algebra $B$ by means of infinitesimal
automorphisms (derivations). This means that we have linear operators $%
\delta _{X}:B\rightarrow B$ which linearly depend on $X$ and satisfy%
\begin{equation*}
\delta _{X}(a\cdot b)=(\delta _{X}a)\cdot b+a\cdot (\delta _{X}b)
\end{equation*}%
for any $a,b\in B.$ The mapping $X\rightarrow \delta _{X}$ is a Lie algebra
homomorphism, i. e. $\delta _{\lbrack X,Y]}=[\delta _{X},\delta _{Y}].$

Now we consider respectively instead of commutative spaces 
 $\mathbf{V}$  and $\beta $
the finite \ projective $\mathcal{A}$--module $\mathcal{E}_{N},$ provided
with N--connection structure, and the finite projective $\mathcal{B}$%
--module $\mathcal{E}_{\beta }.$

A d--connection $\nabla _{X}$ on $\mathcal{E}_{\beta }$ is
by definition a set of linear d--operators, adapted to the N--connection
structure, depending linearly on $X$ and satisfying the Leibniz rule%
\begin{equation}
\nabla _{X}(b\cdot e)=b\cdot \nabla  _{X}(e)+\delta
_{X}b\cdot e  \label{rul1n}
\end{equation}%
for any $e\in \mathcal{E}_{\beta }$ and $b\in \mathcal{B}.$ The rule 
(\ref{rul1n}) is a noncommutative generalization of (\ref{rul1c}). We emphasize
that both operators $\nabla _{X}$ and $\delta _{X}$ are
distinguished by the N--connection structure and that the difference of two
such linear d--operators, $\nabla _{X}-\nabla _{X}^{\prime }$ 
commutes with action of $B$ on $\mathcal{E}_{\beta },$ which
is an endomorphism of $\mathcal{E}_{\beta }.$ Hence, if we fix some fiducial
connection $\nabla _{X}^{\prime }$ (for instance, $%
\nabla _{X}^{\prime }=D_{X})$ on $\mathcal{E}_{\beta }$ an
arbitrary connection has the form 
\begin{equation*}
\nabla _{X}=D_{X}+B_{X},
\end{equation*}%
where $B_{X}\in End_{B}\mathcal{E}_{\beta }$ depend linearly on $X.$

The curvature of connection $\nabla _{X}$ is a two--form $F_{XY}$
which values linear operator in $\mathcal{B}$ and measures a deviation of
mapping $X\rightarrow \nabla _{X}$ from being a Lie algebra
homomorphism,%
\begin{equation*}
F_{XY}=[\nabla _{X},\nabla _{Y}]-\nabla _{\lbrack X,Y]}.
\end{equation*}%
The usual curvature d--tensor is defined as 
\begin{equation*}
F_{\alpha \beta }=\left[ \nabla _{\alpha },\nabla
_{\beta }\right] -\nabla _{\lbrack \alpha ,\beta ]}.
\end{equation*}

The simplest connection on a finite projective $\mathcal{B}$--module $%
\mathcal{E}_{\beta }$ is to be specified by a projector $P:\mathcal{B}%
^{k}\otimes \mathcal{B}^{k}$ when the d--operator $\delta _{X}$ acts
naturally on the free module $\mathcal{B}^{k}.$ The operator $%
\nabla _{X}^{LC}=P\cdot \delta _{X}\cdot P$ $\ $\ is called the
Levi--Civita operator and satisfy the condition
 $Tr[\nabla _{X}^{LC},\phi ]=0$ for any endomorphism
 $\phi \in End_{B}\mathcal{E}_{\beta}.$ 
From this identity, and from the fact that any two connections differ by
an endomorphism that $Tr[\nabla _{X},\phi ]=0$
for an arbitrary connection $\nabla _{X}$ and an arbitrary
endomorphism $\phi ,$ that instead of $\nabla _{X}^{LC}$ we may
consider equivalently the canonical d--connection, constructed only from
d-metric and N--connection coefficients.

\section{Nonholonomic Clifford--Lagrange Structu\-res}

The geometry of spinors on generalized Lagrange and Finsler spaces was
elaborated in Refs. \cite{vsp1,vsph,vst,vsp2}. It was applied for
definition of noncommutative extensions of the Finsler geometry
related to certain models of Einstein, gauge and string gravity
\cite{vstr,vncg,vncgs,vex1,vt,vp}. Recently, it is was proposed an
extended Clifford approach to relativity, strings and noncommutativity
based on the concept of "C--space'' \cite{castro1,castro2,castro3,castro4}.

 The aim of this section is to
 formulate the geometry of nonholonomic Clifford--Lagrange structures
 in a form adapted to generalizations for noncommutative spaces.

\subsection{Clifford d--module} \label{sscldm}
Let  $\mathbf{V}$ be a compact N--anholonomic manifold. We denote,
respectively, by  $T_x\mathbf{V}$ and $T^*_x\mathbf{V}$  the tangent
and cotangent spaces in a point $x\in\mathbf{V}.$ We  consider a 
 complex vector bundle $\tau: E \rightarrow \mathbf{V}$
where, in general, both the base $\mathbf{V}$ and the total space $E$ may be
provided with N--connection structure, and denote by 
$\Gamma ^{\infty}(E)$ (respectively, $\Gamma(E))$ the set of
differentiable (continous) sections of $E.$ The symbols  
$\chi (\mathbf{M})= \Gamma ^{\infty}(\mathbf{TM})$ and 
 $\Omega ^1 (\mathbf{M}) \doteq \Gamma ^{\infty}\mathbf{(T^*M})$ are used
 respectively for the set of d--vectors and one d--forms on 
$\mathbf{TM}.$ 

\subsubsection{Clifford--Lagrange  functionals} 
In the simplest case, a generic nonholonomic Clifford structure
 can be associated  to a Lagrange metric on a $n$--dimensional real
  vector space $V^n$  provided with a Lagrange quadratic form
 $L(y)=q_L(y,y),$ see subsection \ref{sslfm}. We consider the exterior
 algebra $\wedge V^n$ defined by the identity element $\I$ and
 antisymmetric products $v_{[1]}\wedge ... \wedge v_{[k]}$ with 
 $v_{[1]}, ...,  v_{[k]} \in V^n$ for $k \leq dim V^n$ where 
$\I \wedge v = v,$\ $v_{[1]}\wedge v_{[2]}=-v_{[2]}\wedge v_{[1]}, ...$
\begin{definition}\label{dclalg} The Clifford--Lagrange (or Clifford--Minkowski)
 algebra is a $\wedge V^n$ algebra provided with a product 
 \begin{equation}
 uv+vu = 2 ^{(L)}g(u,v)\ \I \label{clpl} 
 \end{equation}
\begin{equation}
 \mbox{(or }uv+vu = 2 ^{(F)}g(u,v)\ \I \ ) \label{clpm}
\end{equation}
for any $u,v \in V^n$ and $\  ^{(L)}g(u,v)$ (or $\  ^{(F)}g(u,v))$
defined by formulas (\ref{lagm}) (or(\ref{finm})). 
\end{definition}

For simplicity, hereafter we shall prefer to  write down  the formulas
for the Lagrange configurations instead of   dubbing of similar
formulas  for the Finsler configurations. 

We can introduce the complex Clifford--Lagrange algebra $\C
l_{(L)}(V^n)$ structure by using the complex unity ``i'',
 $V_{\C} \doteq V^n + i V^n,$ enabled with complex metric
\begin{equation*}
   ^{(L)}g_{\C}(u,v+iw) \doteq \    ^{(L)}g(u,v) + i\    ^{(L)} g(u,w),
\end{equation*}
which results in certain isomorphisms of matix algebras (see, for
instance, \cite{bondia}),
\begin{eqnarray*} 
\C l (\R ^{2m} )&\simeq & M_{2^m}(\C ), \\
 \C l (\R ^{2m+1}) &\simeq & M_{2^m}(\C ) \oplus  M_{2^m}(\C ) .
\end{eqnarray*} 
We omitted the label $(L)$ because such isomorphisms hold true for any
quadratic forms.

The Clifford--Lagrange algebra possesses usual properties:
\begin{enumerate}
\item On $\C l_{(L)}(V^n)$, it is linearly defined the involution "*",
\begin{equation*} {(\lambda v_{[1]} ...  v_{[k]})}^{*} =
 \overline{\lambda} v_{[1]} ...  v_{[k]},\ \forall \lambda \in \C . 
\end{equation*}
\item There is a $\Z _2$ graduation,
\begin{equation*}
\C l_{(L)}(V^n)= \C l_{(L)}^+(V^n) \oplus \C l_{(L)}^{-}(V^n)
\end{equation*}
with $\chi _{(L)}(a) = \pm 1$ for $a \in  \C l_{(L)}^{\pm}(V^n),$ where 
$\C l_{(L)}^+(V^n),$ or  $\C l_{(L)}^{-}(V^n),$ are defined by
products of an odd, or even, number of vectors.
\item For positive definite forms $q_L(u,v)$, one defines the
  chirality of\\  $\C l_{(L)}(V^n),$ 
\begin{equation*}
\gamma _{(L)}= (-i)^n´ e_1 e_2... e_n,\ \gamma ^2 =\gamma ^* \gamma = \I
\end{equation*}
where $\lbrace e_i \rbrace ^n_{i=1}$ is an orthonormal basis of $V^n$
and $n=2n',$ or $=2n'+1.$
\end{enumerate}

In a more general case, a  nonholonomic Clifford structure
  is defined by quadratic d--metric form
  $\mathbf{q}(\mathbf{x},\mathbf{y})$ (\ref{qdmf}) on 
 a $n+m$--dimensional real d--vector space $V^{n+m}$  with the
  $(n+m)$--splitting defined by the N--connection structure.
\begin{definition} \label{dcdalg} The Clifford d--algebra  is a
 $\wedge V^{n+m}$ algebra provided with a product 
 \begin{equation}
 \mathbf{u}\mathbf{v} + \mathbf{v} \mathbf{u} =
 2 \mathbf{g} (\mathbf{u}, \mathbf{v})\ \I \label{cdalg} 
 \end{equation}
or, equivalently, distinguished into h-- and v--products
\begin{equation*}
 uv+vu = 2 g(u,v)\ \I \  
\end{equation*}
and
\begin{equation*}
 \ ^\star u\ ^\star v + \ ^\star v \ ^\star u 
= 2\ ^\star h(\ ^\star u,\ ^\star v)\ \I \  
\end{equation*}
for any $\mathbf{u}=(u,\ ^\star u),\ \mathbf{v} = (v,\ ^\star v) \in
 V^{n+m}.$
\end{definition}
Such Clifford d--algebras have similar properties on the irreducible
h-- and v--components  as the
Clifford--Lagrange algebras. We may define a standard complexification but it
should be emphasized that for $n=m$ the N--connection (in
particular, the canonical Lagrange N--connection) induces naturally
an almost complex structure (\ref{acs1}) which gives the possibility
to define almost complex Clifford d--algebras
 (see details in \cite{vsp1,vsp2}).

\subsubsection{Clifford--Lagrange  and Clifford  N--anholonomic
  structures}

A metric  on a manifold $M$ is defined by sections of the
tangent bundle $TM$ provided with a bilinear symmetric form on
continous sections $\Gamma (TM).$  In Lagrange geometry, the metric
structure is of type $ ^{(L)}g_{ij}(x,y)$ (\ref{lagm}) which allows us to
define Clifford--Lagrange algebras $\C l_{(L)}(T_x M),$ in any point
$x\in TM,$
\begin{equation*}
\gamma _i \gamma _j + \gamma _j \gamma _i = 2\  ^{(L)}g_{ij}\ \I.
\end{equation*}  
For any point $x\in M$ and fixed $y=y_0,$
 one exists a standard complexification,
 $T_x M^{\C} \doteq T_x M + i T_x M,$ which can be used for definition
 of the 'involution' operator on sections of $T_x M^{\C},$
\begin{equation*}
 \sigma _1 \sigma _2 (x) \doteq \sigma _2 (x) \sigma _1 (x), \ 
 \sigma ^* (x) \doteq \sigma (x) ^*, \forall x \in M,
\end{equation*} 
where  "*" denotes the involution on every $\C l_{(L)}(T_x M).$ The
norm is defined by using the Lagrange norm, see Definition
\ref{dlagf}, 
\begin{equation*}
 \parallel \sigma \parallel _L \doteq sup_{x\in M}\ \lbrace 
\mid \sigma (x) \mid _L \rbrace,
\end{equation*}  
which defines a $C^*_{L}$--algebra instead of the usual $C^*$--algebra 
 of $\C l(T_x M).$ Such constructions can be also performed on the
 cotangent space $T_x M,$  or for any  vector bundle $E$ on $M$
 enabled with a symmetric bilinear form of class $C^\infty$ on  
 $\Gamma ^\infty (E) \times \Gamma ^\infty (E).$ 

For Lagrange spaces modeled on $\widetilde{TM},$ there is a natural
almost complex structure $\mathbf{F}$\ (\ref{acs1}) induced by the
canonical N--connection $  ^{(L)}N,$  see the Results \ref{r2}, \ref{r4} and
\ref{r5}, which allows also to construct an almost Kahler model of
Lagrange geometry, see details in Refs. \cite{ma1,ma2}, and to define
an Clifford--Kahler d--algebra $\C l_{(KL)}(T_x M)$ \cite{vsp1}, 
for $y=y_0,$ being provided with the norm
\begin{equation*}
\parallel \sigma \parallel _{KL} \doteq sup_{x\in M}\ \lbrace 
\mid \sigma (x) \mid _{KL} \rbrace,
\end{equation*}   
which on $T_x M$ is defined by projecting on $x$ the d--metric
 $\ ^{(L)}\mathbf{g}$ (\ref{slm}). 

In order to model Clifford--Lagrange structures on $\widetilde{TM}$
and $\widetilde{T^*M}$ it is necessary to consider d--metrics 
induced by Lagrangians:

\begin{definition} A Clifford--Lagrange space on a manifold $M$
  enabled with  a fundamental metric $\ ^{(L)}g_{ij}(x,y)$
 (\ref{lqf})  and canonical  N--connection $\ ^{(L)}N_{j}^{i}$
 (\ref{cncl}) inducing a Sasaki type d--metric $\ ^{(L)}\mathbf{g}$ 
(\ref{slm}) is  defined as a Clifford bundle 
$\C l_{(L)}(M)\doteq \C l_{(L)}(T^*M).$  
\end{definition}

For a general N--anholonomic manifold $\mathbf{V}$ of dimension $n+m$ 
provided with a general  d--metric structure $\mathbf{g}$ (\ref{metr})
 (for instance, in a gravitational model, or constructed by conformal
 transforms and  embeddings into higher dimensions of a Lagrange (or Finsler)
 d--metrics), we introduce 
\begin{definition} A Clifford N--anholonomic bundle on $\mathbf{V}$ is
  defined as\\ $\C l_{(N)}(\mathbf{V})\doteq \C l_{(N)}(T^*\mathbf{V}).$  
\end{definition}

Let us consider a complex vector bundle $\pi : \mathbf{E} \to M$
provided with N--connection structure which can be defined by a
corresponding exact chain of subbundles, or nonintegrable
distributions, like for real vector bundles, see \cite{ma1,ma2} and 
 subsection \ref{snc}.   Denoting by
$V^m_{\C}$ the typical fiber (a complex vector space), we can define
the usual Clifford map
\begin{equation*}
c:\   \C l (T^*M) \rightarrow End(V^m_{\C}) 
\end{equation*} 
via (by convention, left) action on sections
 $c(\sigma) \sigma ^1(x) \doteq c(\sigma (x)) \sigma ^1 (x).$
\begin{definition} \label{defcdma} The Clifford d--module (distinguished by a
  N--connection)  of a N--anholonomic vector bundle $\mathbf{E}$ is
  defined by the $C(M)$--module $\Gamma (\mathbf{E})$ of continuous
  sections in $\mathbf{E},$
 \begin{equation*}
 c:\ \Gamma(\C l (M)) \rightarrow End(\Gamma (\mathbf{E})). 
 \end{equation*}
\end{definition}

In an alternative case,  one considers a complex
 vector  bundle $\pi:\ E \to  \mathbf{V}$   on an  
N--anholonomic space $\mathbf{V}$ when the N--connection structure is
 given for the base manifold.
\begin{definition} \label{defcdmb} The Clifford d--module 
  of a  vector bundle ${E}$ is
  defined by the $C(\mathbf{V})$--module $\Gamma ({E})$ of continuous
  sections in ${E},$  
 \begin{equation*}
 c:\ \Gamma(\C l_{(N)} (\mathbf{V})) \rightarrow End(\Gamma ({E})). 
 \end{equation*}
\end{definition}

A Clifford d--module with both N--anholonomic total space
$\mathbf{E}$ and base space $\mathbf{V}$ with corresponding
N--connections (in general, two independent ones, but the
N--connection in the distinguished complex vector bundle must be
adapted to the N--connection on the base) has to be defined by an
"interference" of Definitions \ref{defcdma} and \ref{defcdmb}.

\subsection{N--anholonomic spin structures}
Usually, the spinor bundle on a manifold $M,\ dim M=n,$ is constructed
on the tangent bundle by substituting the group $SO(n)$ by its
universal covering $Spin(n).$ If a Lagrange fundamental quadratic form  
$\ ^{(L)}g_{ij}(x,y)$ (\ref{lqf}) is defined on $T_x,M$  we can
consider  
Lagrange--spinor spaces in every point $x \in M.$  The constructions
can be completed on $\widetilde{TM}$ by using the Sasaki type metric 
$\ ^{(L)}\mathbf{g}$ (\ref{slm}) being similar for any type of
N--connection and d--metric structure on $TM.$ On general
N--anholonomic manifolds $\mathbf{V}, dim \mathbf{V}=n+m,$  the
distinguished  spinor space (in brief, d--spinor
space) is to be  derived from the d--metric (\ref{metr}) and adapted
to the N--connection structure. In this case, the group $SO(n+m)$ is
not only substituted by $Spin(n+m)$ but with respect to N--adapted
frames (\ref{dder}) and (\ref{ddif}) one defines irreducible
decompositions to $Spin(n)\oplus Spin(m).$

\subsubsection{Lagrange spin groups}

Let us consider a vector space $V^n$ provided with Clifford--Lagrange
structures as in subsection \ref{sscldm}. We denote a such space as
$V^n_{(L)}$ in order to emphasize that  its tangent space 
 is provided with a Lagrange
type quadratic form $\ ^{(L)}g.$ In a similar form, we shall write
 $\C l_{(L)} (V^n)  \equiv \C l (V^n_{(L)})$ if this will be more
 convenient.  A vector $u \in V^n_{(L)}$ has a unity length 
 on the Lagrange quadratic form if $\ ^{(L)}g(u,u) = 1,$ or
 $u^2=\I,$ as an element of corresponding Clifford algebra, 
 which follows from (\ref{clpl}). We define an endomorphism of $V^n:$ 
\begin{equation*}
\phi _{(L)}^{u}\doteq \chi _{(L)} (u) v u^{-1}=-uvu = \left( uv - 2\
^{(L)}g(u,v)\right) u = u - 2\ ^{(L)}g(u,v)u 
\end{equation*}
where $\chi _{(L)}$ is the $\Z _2$ graduation defined by $\ ^{(L)}g.$
 By multiplication, 
\begin{equation*} 
\phi _{(L)}^{u_1 u_2}(v)\doteq u_2^{-1}  u_1^{-1} v u_1 u_2 =
 \phi _{(L)}^{u_2} \circ \phi _{(L)}^{u_1} (v),     
\end{equation*}
which defines the subgroup $SO(V^n_{(L)}) \subset O(V^n_{(L)}).$ Now
we can define \cite{vsp1,vsp2}
\begin{definition} The space of complex Lagrange spins is defined by the 
 subgroup   $Spin^c_{(L)}(n)\equiv Spin^c (V^n_{(L)}) 
 \subset \C l (V^n_{(L)}),$ 
 determined  by the products of pairs of 
  vectors $w \in V_{(L)}^{\C}$ when $w\doteq \lambda u$ where
  $\lambda$ is a complex number of module 1 and $u$ is of unity length
  in $V^n_{(L)}.$
\end{definition}

We note that $ker \phi _{(L)} \cong U(1).$ We can define a
homomorphism $\nu _{(L)}$ with values in $U(1),$
\begin{equation*} 
\nu _{(L)}(w)= w_{2k}...w_1 w_1...w_{2k} = \lambda _1 ... \lambda _{2k},
\end{equation*}    
where $w =  w_1...w_{2k} \in Spin^c (V^n_{(L)})$  
and $\lambda _i = w_i^2 \in U(1).$
\begin{definition}
The group of real Lagrange spins  $Spin^c_{(L)}(n) \equiv Spin
(V^n_{(L)})$ is defined by  $ker\  \nu _{(L)}.$
\end{definition}
The complex conjugation on $ \C l (V^n_{(L)})$  is usually defined as
$\overline{\lambda v} \doteq \overline{\lambda} v$ for $\lambda \in
\C,\ v\in V^n_{(L)}.$ So, any element $w \in Spin (V^n_{(L)})$
satisfies the conditions $\overline{w^*}w = w^* w = \I$  and
$\overline{w}=w.$ If we take $V^n_{(L)}=\R ^n$ provided with a (pseudo)
Euclidean quadratic form instead of the Lagrange norm, we obtain the
usual spin--group constructions from the (pseudo) Euclidean geometry.

\subsubsection{Lagrange spinors and d--spinors:\ Main Result 1}

A usual spinor is a section of a vector bundle $S$ on a manifold $M$
when  an irreducible representation of the group $Spin (M)
\doteq Spin (T^*_x M)$ is defined on the typical fiber. The set of
sections $\Gamma (S)$ is a irreducible Clifford module. If the base
manifold of type $M_{(L)},$ or  is a general N--anholonomic manifold 
$\mathbf{V},$ we have to define the spinors on such spaces as to be
adapted to the respective N--connection structure.

In the case when the base space is of even dimension (the geometric
constructions in in this subsection will be considered for even
dimensions both for the base and typical fiber spaces), one should consider the
so--called Morita equivalence (see details in \cite{bondia,mart2}
 for a such equivalence between $C(M)$ and $\Gamma (\C l(M))$). One
 says that two algebras $\mathcal{A}$ and $\mathcal{B}$ are
 Morita--equivalent if  
\begin{equation*}
\mathcal{E} \otimes _\mathcal{A} \mathcal{F} \simeq \mathcal{B}
 \mbox{ and }
\mathcal{F} \otimes _\mathcal{B} \mathcal{F} \simeq \mathcal{A}, 
\end{equation*}
respectively, for  $\mathcal{B}$-- and $\mathcal{A}$--bimodules and 
 $\mathcal{B-A}$--bimodule $\mathcal{E}$ and
 $\mathcal{A-B}$--bimodule $\mathcal{F}.$
If we study algebras through theirs representations, we also have to
 consider various algebras related by the Morita equivalence.
\begin{definition}
  A Lagrange spinor  bundle $S_{(L)}$ on a manifold $M,\ dim
 M=n,$ is a complex vector bundle with both  defined action of the spin group 
$Spin (V^n_{(L)})$ on the typical fiber  and an irreducible
 reprezentation of the group 
$ Spin_{(L)} (M) \equiv Spin (M_{(L)}) \doteq Spin (T_x^*
 M_{(L)}).$ The set of sections $\Gamma (S_{(L)})$ defines an irreducible 
 Clifford--Lagrange module.
\end{definition}

The so--called "d--spinors" have been introduced for the spaces
provided with N--connection structure \cite{vsp1,vsph,vstr1}:
\begin{definition} \label{ddsp}
 A distinguished spinor (d--spinor) bundle $\mathbf{S}\doteq (S,\
  ^{\star} S)$ on an
  N--anho\-lo\-nom\-ic manifold $\mathbf{V},$ $\ dim \mathbf{V}=n+m,$ is
  a complex vector bundle with a defined action of the spin
  d--group  $Spin\ \mathbf{V} \doteq Spin(V^n) \oplus  Spin(V^m) $ with
  the splitting adapted to the N--connection structure which results
  in an irreducible  representation  
 $Spin (\mathbf{V}) \doteq Spin (T^*\mathbf{V}).$ The set of sections
 $\Gamma (\mathbf{S}) = \Gamma\ ({S}) \oplus 
 \Gamma (\ ^{\star} {S}) $ is an irreducible  Clifford d--module.   
\end{definition}

The fact that $C(\mathbf{V})$ and $\Gamma
(\C l (\mathbf{V}))$ are Morita equivalent can be analyzed by applying
in N--adapted form, both on the base and fiber spaces, the
consequences of the Plymen's theorem (see Theorem 9.3 in
Ref. \cite{bondia}). This is connected with the possibility to
distinguish the $Spin (n)$ (or, correspondingly $Spin (M_{(L)}),$\ 
 $Spin(V^n) \oplus  Spin(V^m) )$ an antilinear bijection $J:\ S\ \to \ S$
 (or $J:\ S_{(L)}\ \to \ S_{(L)}$ and $J:\ \mathbf{S}\ \to \
 \mathbf{S})$ with the properties:
\begin{eqnarray}
J(\psi f) &=& (J\psi )f \mbox{ for } f \in  C(M) (\mbox{ or }
C(M_{(L)}),\  C(\mathbf{V})); \nonumber \\
J(a\psi ) &=& \chi (a)J\psi, \mbox{ for } a\in \Gamma ^\infty (\C
l(M))  (\mbox{ or } \Gamma ^\infty (\C l(M_{(L)})),\  
 \Gamma ^\infty (\C l(\mathbf{V}));   \nonumber \\
(J\phi | J\psi ) &=& (\psi | \phi) \mbox{ for } \phi , \psi \in S 
(\mbox{ or } S_{(L)}, \mathbf{S} ).
 \label{jeq} 
\end{eqnarray}

\begin{definition}
The spin structure on a manifold $M$ (respectively, on 
 $M_{(L)},$ or on N--anholonomic manifold $\mathbf{V}$) with even dimensions
for the corresponding base and typical fiber spaces is defined by a
bimodule $S$ (respectively, $M_{(L)},$ or $\mathbf{V}$) obeying the
Morita equivalence $C(M)- \Gamma (\C l(M))$ (respectively,
 $C(M_{(L)})- \Gamma (\C l(M_{(L)})),$ or 
$C(\mathbf{V})- \Gamma (\C l(\mathbf{V})))$   by a corresponding
bijections (\ref{jeq}) and a fixed orientation on $M$ (respectively,
on $M_{(L)}$ or $\mathbf{V}$). 
\end{definition}

In brief, we may call $M$ ($M_{(L)},$ or $\mathbf{V}$) as a spin manifold
(Lagrange spin manifold, or N--anholonomic spin manifold). If any of
the base or typical fiber spaces is of odd dimension, we may perform
similar constructions by considering $\C l^+$ instead of $\C l.$ 

The considerations presented in this Section consists the proof of the
first main Result of this paper (let us conventionally say that it is
the 7th one after the Results \ref{r1}--\ref{r6}:\
\begin{theorem} \label{mr1} {\bf (Main Result 1)} 
Any regular Lagrangian and/or 
N--con\-nec\-ti\-on structure define naturally the fundamental geometric objects
 and structures  (such as the Clifford--Lagrange module and Clifford d--modules,
 the Lagrange spin structure  and d--spinors) for the corresponding Lagrange
 spin manifold  and/or  N--anholo\-nom\-ic spinor (d--spinor) manifold.
\end{theorem}
We note that similar results were obtained in
Refs. \cite{vsp1,vsph,vst,vsp2} for the standard Finsler and Lagrange
geometries and theirs higher order generalizations. In a more
restricted form,  the idea  of  Theorem \ref{mr1}
can be found in Ref. \cite{vncg}, where the first models of noncommutative
Finsler geometry and related gravity  were considered (in a more rough
form, for instance, with constructions not reflecting the Morita equivalence).

Finally, in this Section, we can make the
\begin{conclusion} Any regular Lagrange and/or N--connection structure
  (the last one being any admissible N--connection in Lagrange--Finsler geometry
  and their generalizations, or induced by any generic off--diagonal and/ or
   nonholonomic frame structure) define certain, corresponding,
  Clifford--La\-gran\-ge module and/or Clifford d--module and  related 
  Lagrange spinor and/or d--spinor structures.  
\end{conclusion}
It is a bit surprizing
  that a Lagrangian may define  not only the fundamental geometric objects of 
 a  nonholonomic Lagrange space but also the structure of 
 a naturally associated Lagrange spin manifold. The Lagrange mechanics and
  off--diagonal gravitational interactions (in general, with nontrivial
  torsion and nonholonomic constraints) can be completely  
 geometrized on Lagrange spin (N--anholonomic) manifolds.  

\section{The Dirac Operator, Nonholonomy, and \newline Spec\-tral
  Trip\-les} \label{sdonst}
The Dirac operator for a certain class of  (non) commutative Finsler
spaces provided with compatible metric structure  was  introduced in
Ref. \cite{vncg} following previous constructions for  the Dirac equations
on locally anisotropic spaces \cite{vsp1,vsph,vstr1,vst,vsp2}. 
The aim of this Section is to elucidate the possibility of  definition of Dirac
operators for general N--anholonomic manifolds and 
Lagrange--Finsler spaces. It should be noted that such geometric
constructions depend on the type of linear connections which are used
  for the complete definition of the  Dirac operator. They are metric
 compatible and N--adapted if the canonical d--connection 
  is used, see Proposition \ref{pcnc} (we can also use  
 any its deformation which results in a  metric compatible
 d--connection). The constructions  can be  more sophisticate and
 nonmetric (with some  geometric  objects not  completely
 defined  on the  tangent spaces) if the Chern, or the 
  Berwald d--connection, is considered, see  Example \ref{ecrdc}. 

\subsection{N--anholonomic Dirac operators}
We introduce the basic definitions and formulas with respect to
N--adapted frames of type (\ref{dder}) and (\ref{ddif}).
 Then we shall present the main  results in a global form.

\subsubsection{Noholonomic vielbeins and spin d--connections}
Let us consider a Hilbert space of finite dimension. For a local dual
coordinate basis $e^{\underline i} \doteq dx^{\underline i}$
 on a manifold  $M,\ dim M =n,$\
  we may respectively introduce certain classes of
    orthonormalized  vielbeins and the
  N--adapted vielbeins, \footnote{(depending both on the base coordinates $x
  \doteq\ x^i$   and some "fiber" coordinates $y\doteq y^a,$ the
  status of $y^a$ depends on what kind of models we shall consider:\  
  elongated on $TM,$ for a Lagrange space, for a vector bundle, or on
  a N--anholonomic manifold)}
\begin{equation}
  e^{\hat i} \doteq  e^{\hat i}_{\ \underline i} (x,y)\ e^{\underline i}
 \mbox{ and }
   e^{i} \doteq e^{i}_{\ \underline i}(x,y)\ e^{\underline i}, \label{hatbvb} 
\end{equation}
where 
\begin{equation*}
 g^{\underline i \underline j}(x,y)\ e^{\hat i}_{\ \underline i} (x,y)
 e^{\hat j}_{\ \underline j} (x,y)
 = \delta ^{\hat i \hat j} \mbox{ and } 
 g^{\underline i \underline j}(x,y)\ e^{i}_{\ \underline i} (x,y)
 e^{j}_{\ \underline j} (x,y)
 = g ^{ij}(x,y).
\end{equation*}
We define the the algebra of Dirac's gamma matrices (in brief,
h--gamma matrices defined by  self--adjoints
matrices $M_k(\C)$ where $k=2^{n/2}$ is the dimension of the
irreducible representation of $\C l(M)$ for even dimensions, or of $\C
l(M)^+$ for odd dimensions) from the relation
\begin{equation} 
 \gamma ^{\hat i} \gamma ^{\hat j} +  \gamma ^{\hat j}\gamma ^{\hat i}
 = 2 \delta ^{\hat i \hat j}\ \I. \label{grelflat}
\end{equation}
We can consider the action of $dx^i\in \C l (M)$ on a spinor $\psi \in
S$ via representations 
\begin{equation}
\ ^{-}c(dx^{\hat i}) \doteq \gamma ^{\hat i} \mbox{ and }
 \ ^{-}c(dx^i)\psi \doteq \gamma ^i
 \psi \equiv e^i_{\ \hat i}\ \gamma ^{\hat i} \psi. \label{gamfibb}
\end{equation}

For any type of spaces $T_xM, TM, \mathbf{V}$ possessing a local (in
any point) or global fibered structure and, in general, enabled with
a N--connection structure, we can introduce similar definitions of the
gamma matrices following algebraic relations and metric structures on
fiber subspaces,
\begin{equation}
  e^{\hat a} \doteq  e^{\hat a}_{\ \underline a} (x,y)\ e^{\underline a}
 \mbox{ and }
   e^{a} \doteq e^{a}_{\ \underline a}(x,y)\ e^{\underline a}, \label{hatbvf} 
\end{equation}
where 
\begin{equation*}
 g^{\underline a \underline b}(x,y)\ e^{\hat a}_{\ \underline a} (x,y)
 e^{\hat b}_{\ \underline b} (x,y)
 = \delta ^{\hat a \hat b} \mbox{ and } 
 g^{\underline a \underline b}(x,y)\ e^{a}_{\ \underline a} (x,y)
 e^{b}_{\ \underline b} (x,y)
 = h^{ab}(x,y).
\end{equation*}
Similarly, we define the algebra of Dirac's matrices related to
typical fibers (in brief, v--gamma matrices defined by self--adjoints
matrices $M_k'(\C)$ where $k'=2^{m/2}$ is the dimension of the
irreducible representation of $\C l(F)$ for even dimensions, or of $\C
l(F)^+$ for odd dimensions, of the typical fiber) from the relation
\begin{equation} 
 \gamma ^{\hat a} \gamma ^{\hat b} +  \gamma ^{\hat b}\gamma ^{\hat a}
 = 2 \delta ^{\hat a \hat b}\ \I. \label{grelflatf}
\end{equation}
The action of $dy^a\in \C l (F)$ on a spinor $\ ^\star
\psi \in \ ^{\star}S$ is considered via representations 
\begin{equation}
\ ^{\star}c(dy^{\hat a}) \doteq \gamma ^{\hat a} \mbox{ and }
 \ ^{\star}c(dy^a)\ ^\star \psi \doteq \gamma ^a
 \ ^\star \psi \equiv e^a_{\ \hat a}\ \gamma ^{\hat a} 
 \ ^\star \psi. \label{gamfibf}
\end{equation}
 
We note that additionally to formulas (\ref{gamfibb}) and
(\ref{gamfibf}) we may write respectively 
\begin{equation*}
 c(dx^{\underline i})\psi \doteq \gamma ^{\underline i}
 \psi \equiv e^{\underline i}_{\ \hat i}\ \gamma ^{\hat i} \psi
\mbox{ and }
 c(dy^{\underline a})\ ^\star \psi \doteq \gamma ^{\underline a} 
 \ ^\star \psi \equiv e^{\underline a}_{\ \hat a}\ \gamma ^{\hat a} 
 \ ^\star \psi
\end{equation*}
but such operators are not adapted to the N--connection structure.

A more general gamma matrix calculus with distinguished gamma
matrices (in brief, d--gamma matrices\footnote{in our previous works
  \cite{vsp1,vsph,vstr1,vst,vsp2} we wrote $\sigma$ instead of
  $\gamma$}) can be elaborated for N--anholonomic manifolds
$\mathbf{V}$ provided with d--metric structure $\mathbf{g}=[g,
^{\star}g]$  and for
 d--spinors $\breve{\psi} \doteq (\psi,\ ^{\star}\psi)
\in \mathbf{S} \doteq (S,\ ^{\star}S),$ see the corresponding 
 Definitions \ref{dnam}, \ref{ddms} and \ref{ddsp}. Firstly, we should write
 in a unified form, related to a d--metric (\ref{metr}), the formulas
 (\ref{hatbvb}) and (\ref{hatbvf}),  
 \begin{equation}
  e^{\hat \alpha} \doteq  e^{\hat \alpha}_{\ \underline a}
 (u)\ e^{\underline \alpha}
 \mbox{ and }
   e^{\alpha} \doteq e^{\alpha}_{\ \underline \alpha}(u)\ 
 e^{\underline \alpha}, \label{hatbvd} 
\end{equation}
where 
\begin{equation*}
 g^{\underline{\alpha} \underline{\beta}}
(u)\ e^{\hat \alpha}_{\ \underline \alpha} (u)
 e^{\hat \beta}_{\ \underline \beta} (u)
 = \delta ^{\hat \alpha \hat \beta} \mbox{ and } 
 g^{\underline \alpha \underline \beta}(u)\ e^{\alpha}_{\ \underline
   \alpha}  (u)  e^{\beta}_{\ \underline \beta} (u)
 = g^{\alpha \beta}(u).
\end{equation*}
The second step, is to consider d--gamma matrix relations (unifying
(\ref{grelflat}) and (\ref{grelflatf})) 
\begin{equation} 
 \gamma ^{\hat \alpha} \gamma ^{\hat \beta} +  \gamma ^{\hat \beta}
\gamma ^{\hat \alpha}
 = 2 \delta ^{\hat \alpha \hat \beta}\ \I,  \label{grelflatd}
\end{equation}
 with the action of $du^\alpha \in \C l (\mathbf{V})$ on a d--spinor
 $\breve{\psi} \in \ \mathbf{S}$ resulting in distinguished irreducible
representations (unifying (\ref{gamfibb}) and (\ref{gamfibf}))  
\begin{equation}
\mathbf{c}(du^{\hat \alpha}) \doteq \gamma ^{\hat \alpha} \mbox{ and }
 \mathbf{c}=(du^\alpha)\  \breve{\psi} \doteq \gamma ^\alpha 
 \  \breve{\psi} \equiv e^\alpha_{\ \hat \alpha}\ \gamma ^{\hat \alpha} 
 \  \breve{\psi} \label{gamfibd}
\end{equation}
which allows to write
\begin{equation} 
 \gamma ^{\alpha}(u) \gamma ^{\beta}(u) + 
 \gamma ^{\beta} (u) \gamma ^{\alpha} (u)
 = 2 g ^{\alpha \beta}(u)\ \I.  \label{grelnam}
\end{equation}
In the canonical representation we can write in irreducible form
$\breve \gamma \doteq \gamma \oplus\ ^\star \gamma$ and
 $\breve \psi \doteq \psi \oplus\ ^\star \psi,$  for instance, by
 using block type of h-- and v--matrices, or, writing alternatively
 as couples of gamma and/or h-- and v--spinor objects written 
in N--adapted form,
\begin{equation} \label{crgs}
\gamma ^{\alpha} \doteq (\gamma ^i, \gamma ^a) \mbox{ and } 
 \breve \psi \doteq (\psi,\ ^\star \psi).
\end{equation}
The decomposition (\ref{grelnam}) holds with respect to a N--adapted
vielbein (\ref{dder}). We also note that for a spinor calculus, the
indices of spinor objects should be treated as abstract spinorial ones
 possessing certain reducible, or irreducible, properties depending on
 the space dimension (see details in Refs.
 \cite{vsp1,vsph,vstr1,vst,vsp2}). For simplicity, we shall consider
 that spinors like $\breve \psi, \psi,\ ^\star \psi$ and all type of
 gamma objects can be enabled with corresponding spinor indices
 running certain values which are different from the usual coordinate
 space indices. In a "rough" but brief form we can use the same indices $i,j,
 ...,a,b...,\alpha, \beta,...$ both for d--spinor and d--tensor
 objects.

The spin connection $\nabla ^S$  for the Riemannian manifolds is
induced by the Levi--Civita connection $\ ^\nabla \Gamma ,$
\begin{equation} \label{sclcc}
 \nabla ^S  \doteq d - \frac{1}{4}\ ^\nabla \Gamma ^i_{\ j k} 
\gamma _i \gamma ^j\ dx^k.
\end{equation}
On N--anholonomic spaces, it is possible to define spin connections
which are N--adapted by replacing the Levi--Civita connection by any
d--connection (see Definition \ref{ddc}).
\begin{definition}
The canonical spin d--connection is defined by the canonical
d--connection (\ref{cdc}) as
\begin{equation} \label{csdc}
 \widehat{\nabla}^{\mathbf{S}} \doteq \delta
 - \frac{1}{4}\ \widehat{\mathbf{\Gamma}}^\alpha_{\ \beta \mu} 
\gamma _\alpha \gamma ^\beta \delta u^\mu,
\end{equation}
where the absolute differential $\delta$ acts in  N--adapted form
resulting in 1--forms decomposed with respect to N--elongated
differentials like $\delta u^\mu = (dx^i, \delta y^a)$ (\ref{ddif}).
\end{definition}
We note that the canonical spin d--connection 
$\widehat{\nabla}^{\mathbf{S}}$ is metric compatible and
contains nontrivial d--torsion coefficients
 induced by the N--anholonomy relations (see
the formulas (\ref{dtors}) proved for arbitrary d--connection). 
 It is possible to introduce more general spin d--connections
 ${\mathbf{D}}^{\mathbf{S}}$ by using the same formula (\ref{csdc})
 but for arbitrary metric compatible d--connection 
${\mathbf{\Gamma}}^\alpha_{\ \beta \mu}.$ 

In a particular case, we
 can define, for instance, the canonical spin d--connections for a
local modelling of a $\widetilde{TM}$ space on $\widetilde{\mathbf{V}}%
_{(n,n)}$ with the canonical d--connection 
 $\widehat{\mathbf{\Gamma }}_{\ \alpha \beta
}^{\gamma }=(\widehat{L}_{jk}^{i},\widehat{C}_{jk}^{i}),$ see formulas
 (\ref{candcon1}). This allows us to prove (by considering d--connection
 and d--metric structure defined by the fundamental Lagrange, or Finsler,
 functions, we put formulas (\ref{cncl}) and (\ref{slm})  into
 (\ref{candcon1})):
 \begin{proposition} \label{pcslc} On Lagrange   spaces, there is a
 canonical spin d--connec\-ti\-on (the canonical spin--Lagrange connection),
 \begin{equation} \label{fcslc}
 \widehat{\nabla}^{(SL)} \doteq \delta
 - \frac{1}{4}\ ^{(L)}{\mathbf{\Gamma}}^\alpha_{\ \beta \mu} 
\gamma _\alpha \gamma ^\beta \delta u^\mu,
 \end{equation}
 where $\delta u^\mu = (dx^i, \delta y^k = dy^k +\ ^{(L)} N^k_{\ i}\ dx^i).$
 \end{proposition}
 We emphasize that even regular Lagrangians of classical mechanics
 without spin particles induce in a canonical (but nonholonomic) form
 certain classes of  spin d--connections like (\ref{fcslc}).

For the spaces provided with generic
 off--diagonal metric structure (\ref{ansatz}) (in particular, for
 such  Riemannian manifolds) resulting in  equivalent
 N--anho\-lo\-nom\-ic manifolds, it is possible to 
prove a result being similar to
 Proposition \ref{pcslc}:
\begin{remark} \ 
There is a canonical spin d--connection (\ref{csdc})
 induced by the off--diagonal metric
coefficients with nontrivial $N^a_i$ and associated nonholonomic
 frames in gravity theories.  
\end{remark}

The N--connection structure  also states a global h-- and v--splitting
of spin d--connection operators, for instance, 
\begin{equation} \label{cslc}
 \widehat{\nabla}^{(SL)} \doteq \delta
 - \frac{1}{4}\ ^{(L)}{\widehat{L}}^i_{\ jk } 
\gamma _i \gamma ^j dx^k
 - \frac{1}{4}\ ^{(L)}{\widehat{C}}^a_{\ bc} 
\gamma _a \gamma ^b \delta y^c.
 \end{equation}
So, any spin d--connection is a  d--operator with
conventional splitting of action like ${\nabla}^{(\mathbf{S})} \equiv
 ({\ ^{-}{\nabla}}^{(\mathbf{S})},{\ ^\star
   {\nabla}}^{(\mathbf{S})}),$  or  ${\nabla}^{(SL)} \equiv
 ({\ ^{-}{\nabla}}^{(SL)},{\ ^\star
   {\nabla}}^{(SL)}).$ For instance, for  
$\widehat{\nabla}^{(SL)} \equiv  ({\ ^{-}\widehat{\nabla}}^{(SL)},{\ ^\star
   \widehat{\nabla}}^{(SL)}),$   the operators 
 $\ ^{-}\widehat{\nabla}^{(SL)}$ and 
$\ ^\star \widehat{\nabla}^{(SL)}$ act respectively
  on a h--spinor $\psi$ as 
\begin{equation} \label{hdslop}
{\ ^{-}\widehat{\nabla}}^{(SL)} \psi \doteq  dx^i \  \frac{\delta \psi}{dx^i}
 -  dx^k  \frac{1}{4}\ ^{(L)}{\widehat{L}}^i_{\ jk } 
\gamma _i \gamma ^j \ \psi
\end{equation}
and
\begin{equation*} 
{\ ^\star \widehat{\nabla}}^{(SL)} 
\psi \doteq  \delta y^a \  \frac{\partial \psi}{dy^a}
 - \delta y^c \  \frac{1}{4}\ ^{(L)}{\widehat{C}}^a_{\ bc} 
\gamma _a \gamma ^b \ \psi 
\end{equation*}
being defined by the canonical d--connection (\ref{candcon1}).
\begin{remark} \label{rscfs} We can consider that the h--operator (\ref{hdslop})
  defines a spin generalization of the Chern's d--connection 
 $\ ^{[Chern]}{\mathbf{\Gamma }}_{\ \alpha \beta}^{\gamma }=
(\widehat{L}_{jk}^{i},{C}_{jk}^{i}= 0),\ $ see Example \ref{ecrdc},
 which may be introduced as a minimal extension, with Finsler structure,
 of the spin connection defined by the Levi--Civita connection 
(\ref{sclcc})  preserving the torsionless condition. This is an
 example of nonmetric spin connection operator because 
$\ ^{[Chern]}{\mathbf{\Gamma }}_{\ \alpha \beta}^{\gamma }$ does not
 satisfy the condition of  metric compatibility. 
\end{remark}
We can define  spin Chern--Finsler  structures, considered in
 the Remark \ref{rscfs},  for  any point
 of an N--anholonomic  manifold. There are necessary 
 some additional assumptions in order to completely define such  structures
 (for instance, on the tangent bundle).  We can say that this is
 a  deformed nonholonomic spin
 structure derived from a d--spinor one provided with the canonical
 spin d--connection by deforming the canonical  d--connection
  in a  manner that  the horizontal  torsion vanishes transforming 
 into a nonmetricity d--tensor. The
 "nonspinor" aspects of such generalizations of the Riemann--Finsler
 spaces and gravity models with nontrivial nonmetricity are analyzed
 in Refs. \cite{v2}.
 
\subsubsection{Dirac d--operators:\ Main Result 2}

We consider a vector bundle $\mathbf{E}$ on an N--anholonomic manifold
 $\mathbf{M}$ (with two compatible N--connections defined as h-- and
 v--splittings of $T\mathbf{E}$ and $T\mathbf{M}$)). A d--connection
\begin{equation*} \mathcal{D}:\ \Gamma ^\infty (\mathbf{E})
 \rightarrow \Gamma ^\infty (\mathbf{E}) \otimes \Omega ^1(\mathbf{M})  
\end{equation*}
preserves by parallelism splitting of the tangent total and base
spaces and satisfy the Leibniz condition
\begin{equation*} \mathcal{D}(f \sigma ) = f (\mathcal{D} \sigma) +
  \delta f \otimes \sigma
\end{equation*}
for any $f\in C^\infty (\mathbf{M}),$ and 
 $\sigma \in \Gamma ^\infty (\mathbf{E})$ and $\delta$ defining an
 N--adapted  exterior calculus by using N--elongated operators
 (\ref{dder}) and (\ref{ddif}) which emphasize d--forms instead of
 usual forms on $\mathbf{M},$ with the coefficients taking values in 
$\mathbf{E}.$ 

The metricity and Leibniz conditions for $ \mathcal{D}$ are written
respectively 
\begin{equation} \mathbf{g} (\mathcal{D}\mathbf{X}, \mathbf{Y})+
 \mathbf{g} (\mathbf{X}, \mathcal{D} \mathbf{Y})=
 \delta [\mathbf{g} (\mathbf{X}, \mathbf{Y})], \label{mc1}
\end{equation}
for any $\mathbf{X},\ \mathbf{Y} \in \chi(\mathbf{M}),$
and 
\begin{equation} \mathcal{D} (\sigma \beta) \doteq 
 \mathcal{D} (\sigma) \beta +  \sigma \mathcal{D}(\beta), \label{lc1}
\end{equation}
for any $\sigma, \beta \in \Gamma ^\infty (\mathbf{E}).$

For local computations, we may define the corresponding coefficients of
the geometric d--objects and write 
\begin{equation*} 
\mathcal{D} \sigma _{\acute \beta}
 \doteq {\mathbf{\Gamma}}^{\acute \alpha}_{\ {\acute   \beta} \mu}\  
\sigma _{\acute \alpha} \otimes \delta u^\mu = 
{\mathbf{\Gamma}}^{\acute \alpha}_{\ {\acute   \beta} i}\  
\sigma _{\acute \alpha} \otimes dx^i +  
{\mathbf{\Gamma}}^{\acute \alpha}_{\ {\acute   \beta} a}\  
\sigma _{\acute \alpha} \otimes \delta y^a,
\end{equation*}
where fiber "acute" indices, in their turn, may  split  
${\acute \alpha}\doteq ({\acute i}, {\acute a} )$ if any N--connection
structure is defined on $T\mathbf{E}.$ For some particular
constructions of particular interest, we can take $\mathbf{E} =
T^*\mathbf{V}, = T^*V_{(L)}$ and/or any Clifford d--algebra 
$\mathbf{E} = \C l(\mathbf{V}), \C l(V_{(L)}),...$ with a
corresponding treating of "acute" indices to of d--tensor and/or
d--spinor type as well  when the d--operator  $\mathcal{D}$ transforms into
respective d--connection $\mathbf{D}$ and spin d--connections  
 $\widehat{\nabla}^{\mathbf{S}}$\  (\ref{csdc}),
   $\widehat{\nabla}^{(SL)}$\ (\ref{fcslc}).... All such, adapted to the
   N--connections,  computations are similar for both N--anholonomic 
 (co) vector and spinor  bundles.

The respective actions of the Clifford d--algebra and Clifford--Lagrange algebra
(see Definitions \ref{dcdalg} and \ref{dcdalg}) can be transformed into
 maps  $\Gamma ^\infty (\mathbf{S}) \otimes \Gamma ^(\C
l(\mathbf{V}))$ and $\Gamma ^\infty (S_{(L)}) \otimes \Gamma ^(\C
l(V_{(L)}))$ to $\Gamma ^\infty (\mathbf{S})$ and, respectively, 
 $\Gamma ^\infty (S_{(L)})$ by considering maps of type
 (\ref{gamfibb}) and (\ref{gamfibd})
 \begin{equation*} 
 \widehat{\mathbf{c}}(\breve{\psi} \otimes \mathbf{a}) \doteq
 \mathbf{c}(\mathbf{a}) \breve{\psi}
 \mbox{\ and\ }  
 \widehat{c}({\psi} \otimes {a}) \doteq {c}({a}){\psi}. 
 \end{equation*}
\begin{definition} \label{dddo} The Dirac d--operator (Dirac--Lagrange
 operator) on  a spin N--anholonomic manifold $(\mathbf{V},\mathbf{S}, J)$ (on a
 Lagrange spin manifold\\ $(M_{(L)}, S_{(L)}, J))$  is defined 
\begin{eqnarray}
 \D &\doteq &  -i\ (\widehat{\mathbf{c}} \circ \nabla ^{\mathbf{S}}) 
\label{ddo} \\ 
 & = & \left( \ ^{-}\D =  -i\ (\ ^{-}\widehat{{c}} 
 \circ \ ^{-}\nabla ^{\mathbf{S}}) 
,\ ^\star \D  = -i\ (\ ^{\star}\widehat{{c}} 
 \circ \ ^{\star}\nabla ^{\mathbf{S}} ) \right) \nonumber  \\
 (\ _{(L)}\D  &\doteq & 
 -i\ (\widehat{c} \circ \nabla ^{(SL)}) \ ) \label{dlo} \\
 & = & \left( \ _{(L)}{^{-}\D 
= -i (\ ^{-}\widehat{c} \circ \ ^{-}\nabla ^{(SL)}}),
\ _{(L)}{^\star \D} =  -i (\ ^{\star} \widehat{c} 
 \circ  \ ^{\star}\nabla ^{(SL)}  ) \right)\ ). \nonumber
\end{eqnarray}
Such N--adapted Dirac d--operators are called  canonical  and denoted 
 $\widehat{\D} = ( \ ^{-}\widehat{\D},\
^\star\widehat{\D}\ )$\ ( $ _{(L)}\widehat{\D} = ( \
_{(L)}{^{-}\widehat{\D}},\  _{(L)}{^{\star}\widehat{\D}}\ )$\ )  if they
are defined for the canonical d--connection (\ref{candcon}) (\
(\ref{candcon1})) and respective spin d--connection (\ref{csdc}) 
(\ (\ref{fcslc})).
\end{definition}

Now we can formulate the
 \begin{theorem} \label{mr2} {\bf  (Main Result 2)} Let $(\mathbf{V},\mathbf{S},
   J)$ (\ $(M_{(L)}, S_{(L)}, J)$ be a spin N--anholonomic manifold ( spin
 Lagrange space). There is the canonical Dirac d--operator
 (Dirac--Lagrange operator) defined by the almost Hermitian  spin
  d--operator 
\begin{equation*}
\widehat{\nabla}^{\mathbf{S}}:\ \Gamma ^\infty (\mathbf{S})\rightarrow
 \Gamma ^\infty (\mathbf{S}) \otimes \Omega ^1 (\mathbf{V})
\end{equation*} 
 (spin Lagrange operator  
\begin{equation*}
\widehat{\nabla}^{(SL)}:\ \Gamma ^\infty ({S_{(L)}})\rightarrow
 \Gamma ^\infty (S_{(L)}) \otimes \Omega ^1 (M_{(L)}) \ )
\end{equation*}
 commuting with $J$ (\ref{jeq}) and satisfying the  conditions 
\begin{equation} \label{scmdcond}
 (\widehat{\nabla}^{\mathbf{S}} \breve{\psi}\ |\ \breve{\phi} )\  +
 ( \breve{\psi}\ |\ \widehat{\nabla}^{\mathbf{S}} \breve{\phi} )\ =
 \delta ( \breve{\psi}\ |\  \breve{\phi} )\
\end{equation}
and
\begin{equation*}
\widehat{\nabla}^{\mathbf{S}}( \mathbf{c}(\mathbf{a}) \breve{\psi})\
 = \mathbf{c} (\widehat{\mathbf{D}}\mathbf{a}) \breve{\psi} 
 + \mathbf{c} (\mathbf{a})\widehat{\nabla}^{\mathbf{S}} \breve{\psi}
\end{equation*} 
for $\mathbf{a} \in \C l(\mathbf{V})$ and $\breve{\psi}\in \Gamma
^\infty (\mathbf{S})$ 
\begin{equation} \label{scmlcond}
 (\ (\widehat{\nabla}^{(SL)} \breve{\psi}\ |\ \breve{\phi} )\  +
 ( \breve{\psi}\ |\ \widehat{\nabla}^{(SL)} \breve{\phi} )\ =
 \delta ( \breve{\psi}\ |\  \breve{\phi} )\
\end{equation}
and
\begin{equation*}
\widehat{\nabla}^{(SL)}( \mathbf{c}(\mathbf{a}) \breve{\psi})\
 = \mathbf{c} (\widehat{\mathbf{D}}\mathbf{a}) \breve{\psi} 
 + \mathbf{c} (\mathbf{a})\widehat{\nabla}^{(SL)} \breve{\psi}
\end{equation*} 
for $\mathbf{a} \in \C l(M_{(L)})$ and $\breve{\psi}\in \Gamma
^\infty (S_{(L)}$ )\
 determined by the metricity (\ref{mc1}) and Leibnitz (\ref{lc1})
conditions.
 \end{theorem}
\begin{proof} We sketch the main ideas of  such Proofs. There two
  ways:

 The first one  is similar to
  that given in Ref. \cite{bondia}, Theorem 9.8, for the Levi--Civita
  connection, see similar considerations in \cite{schroeder}. 
In our case, we have to extend the constructions for
  d--metrics and canonical d--connections by applying N--elongated
  operators for differentials and partial derivatives. The
  formulas have to be distinguished into h-- and
  v--irreducible components. We are going to present the related technical details
  in our further publications.

In other turn, the second way, is to argue  a such proof is
 a straightforward consequence of the Result \ref{r6} stating that any
  Riemannian manifold can be modeled as a N--anholonomic manifold
 induced by the generic off--diagonal metric structure. If the results
 from  \cite{bondia} hold true for any Riemannian space,
 the formulas may be rewritten with respect to any local frame
 system, as well with respect to (\ref{dder}) and
 (\ref{ddif}). Nevertheless, on N--anholonomic manifolds the
 canonical d--connection is not just the Levi--Civita connection but
 a deformation of type (\ref{cdc}):\ we must verify that such
 deformations results in  N--adapted constructions satisfying the metricity
 and Leibnitz conditions. The existence of such configurations was
 proven from the properties of the canonical d--connection completely 
 defined from the d--metric and N--connection coefficients. The main
 difference from the case of the Levi--Civita configuration is that we
 have a nontrivial torsion induced by the frame nonholonomy. But it is not
 a problem to define the  Dirac operator with nontrivial torsion if
 the metricity conditions are satisfied. $\Box$
\end{proof}

The canonical Dirac d--operator has very similar
  properties for spin N--anholonomic manifolds and spin Lagrange
  spaces. Nevertheless, theirs geometric and physical meaning may be
  completely different and that why we have written the corresponding
  formulas with different labels and emphasized the existing
  differences. With respect to the {\bf Main Result 2}, one holds 
 three important remarks:
\begin{remark}  The first type of canonical Dirac d--operators may be
  associated to Riemannian--Cartan (in particular, Riemann)
 off--diagonal metric and nonholonomic frame   structures and the
  second type of canonical Dirac--Lagrange operators are completely
  induced by a regular Lagrange mechanics.
  In both cases, such d--operators are
  completely determined by the coefficients of the corresponding
  Sasaki type d--metric and the N--connection structure.
\end{remark}
\begin{remark} The conditions of the Theorem \ref{mr2} may be revised
  for any d--connection and induced spin d--connection satisfying the
  metricity condition. But, for such cases, the corresponding Dirac d--operators are 
  not  completely defined by the d--metric and N--connection
  structures. We can prescribe certain type of torsions of
  d--connections and, via such 'noncanonical' Dirac operators, we are able
  to define noncommutative geometries with prescribed d--torsions. 
\end{remark}
\begin{remark} The properties (\ref{scmdcond}) and  (\ref{scmlcond})
  hold if and only if the metricity conditions are satisfied
  (\ref{mc1}). So, for the Chern or Berwald type d--connections which
  are nonmetric (see Example \ref{ecrdc}  and Remark \ref{rscfs} ),
 the conditions of Theorem \ref{mr2} do not hold. 
\end{remark}
 It is a more   sophisticate problem to find applications in physics for such nonmetric 
  constructions \footnote{See Refs. \cite{hehl} and \cite{v1,v2,v3}
  for details on elaborated geometrical and physical models being, respectively, locally
  isotropic and locally anisotropic.}  but they define positively some examples of nonmetric
  d--spinor and noncommutative structures minimally deformed from the
  Riemannian (non) commutative  geometry to certain Finsler type (non)
  commutative geometries.

\subsection{Distinguished spectral triples}
 The geometric information of a spin manifold (in particular, the
 metric) is contained in the Dirac operator. For nonholonomic
 manifolds, the canonical Dirac d--operator has h-- and v--irreducible parts
 related to  off--diagonal metric terms and nonholonomic frames with
 associated structure. In a more special case, the canonical
 Dirac--Lagrange operator is defined by a regular Lagrangian. So, such
 Driac d--operators contain more information than the usual,
 holonomic, ones.

For simplicity, hereafter, we shall formulate the results for the
general N--anholonomic spaces, by omitting the explicit formulas and
 proofs  for Lagrange and  Finsler spaces, which can be derived by
 imposing certain  conditions  that the 
 N--connection,  d--connection and d--metric are just those
  defined canonically by a  Lagrangian. We shall only present the
  Main Result and some  important Remarks concerning 
 Lagrange   mechanics and/or Finsler  structures.
\begin{proposition} \label{pdohv}
If $\widehat{\D} = ( \ ^{-}\widehat{\D},\ ^\star\widehat{\D}\ )$\ is
the canonical Dirac  d--operator then
\begin{eqnarray*}
\left[\widehat{\D},\ f \right] & = &  i \mathbf{c} (\delta f), 
\mbox{ equivalently,} \\
    \left[\ ^{-}\widehat{\D},\  f \right]  + 
\left[\  ^\star\widehat{\D} ,\ f \right] & = & 
 i \ ^{-}c ( dx^i \frac{\delta  f}{\partial x^i})
 + i \ ^{\star}c( \delta y^a \frac{\partial f}{\partial  y^a}),
\end{eqnarray*}
for all $f \in C^\infty (\mathbf{V}).$
\end{proposition}
\begin{proof} It is a straightforward computation following from
  Definition \ref{dddo}.
\end{proof}

The canonical Dirac d--operator and its irreversible h-- and
v--components have all the properties of the usual Dirac operators
(for instance, they are self--adjoint but unbounded). It is possible
to define a scalar product on $\Gamma ^\infty (\mathbf{S})$,
\begin{equation} \label{scprod}
 <\breve{\psi}, \breve{\phi}> 
\doteq \int_\mathbf{V}(\breve{\psi} | \breve{\phi}) |\nu _\mathbf{g}|
\end{equation}
where $$\nu _\mathbf{g} = \sqrt{det g}\ \sqrt{det h}\ dx^1...dx^n\  
 dy^{n+1}...dy^{n+m}$$ is the volume d--form on the N--anholonomic
 manifold $\mathbf{V}.$

We denote by 
\begin{equation} \label{dhs}
 \mathcal{H}_{N} \doteq L_2 (\mathbf{V}, \mathbf{S}) = \left[
  \ ^{-}\mathcal{H} = L_2  (\mathbf{V},\ ^{-}S),\ ^\star \mathcal{H} =
  L_2 (\mathbf{V},\ ^\star S)\right] 
\end{equation}
the Hilbert d--space obtained by completing $\Gamma ^\infty
(\mathbf{S})$ with the norm defined by the scalar product
(\ref{scprod}).

Similarly to the holonomic spaces, by using formulas (\ref{ddo}) and
 (\ref{csdc}), one may prove that there is a self--adjoint unitary
 endomorphism $\Gamma ^{[cr]}$ of $ \mathcal{H}_{N},$ called
 "chirality", being a ${\Z}_2$ graduation of $ \mathcal{H}_{N},$ 
 \footnote{We use the  label\ $[cr]$\ in order to avoid
 misunderstanding  with the
 symbol $\Gamma$ used for the connections. }
 which satisfies the condition
 \begin{equation}
\widehat{\D}\ \Gamma ^{[cr]} = - \Gamma ^{[cr]}\ \widehat{\D}
\label{chrc}. 
 \end{equation}
We note that the condition (\ref{chrc}) may be written also for 
the irreducible components $\ ^{-}\widehat{\D}$ and $\ ^{\star}\widehat{\D}.$
\begin{definition}
A distinguished canonical spectral triple (canonical spectral d--triple)
 $(\mathcal{A}, \mathcal{H}_{N},\ \widehat{\D}  )$ for an algebra $\mathcal{A}$
 is defined by a Hilbert d--space $\mathcal{H}_{N},$ a representation
 of $\mathcal{A}$ in the algebra $\mathcal{B}(\mathcal{H})$ of d--operators
 bounded on $\mathcal{H}_{N},$ and by a self--adjoint d--operator
 $\widehat{\D},$ of compact resolution,\footnote{An operator $D$ is of
 compact resolution if  for any $\lambda \in sp(D)$ the operator
 $(D-\lambda \I)^{-1}$ is compact, see details in
 \cite{mart2,bondia}.}
 such that
 $[\widehat{\D},a] \in \mathcal{B}(\mathcal{H})$ for any $a\in \mathcal{A}.$ 
\end{definition}

Roughly speaking, every canonical spectral d--triple is defined by two
usual spectral triples which in our case corresponds to certain h--
and v--irreducible components induced by the corresponding h-- and
v--components of the Dirac d--operator. For such spectral h(v)--triples
we can define the notion of $KR^n$--cycle ($KR^m$--cycle) and consider
respective Hochschild complexes.  We note that in order to define a
noncommutative geometry  the h-- and v-- components of a 
 canonical spectral d--triples must satisfy some well defined
 Conditions \cite{connes1,bondia} (Conditions 1 - 7, enumerated in
 \cite{mart2}, section II.4) which states:\ 1) the spectral dimension, being of
 order $1/(n+m)$ for a Dirac d--operator, and of order $1/n$ (or $1/m)$
 for its h-- (or v)--components; 2) regularity; 3) finitness; 4)
 reality; 5) representation of 1st order; 6) orientability; 7)
 Poincar\'{e} duality. Such conditions can be satisfied by any Dirac
 operators and canonical Dirac d--operators (in the last case we have
 to work with d--objects). \footnote{ We omit in this paper the
   details on axiomatics and related proofs for  such
   considerations:\ we shall present details and proofs in our further
   works. Roughly speaking, we are in right to do this because the
   canonical d--connection and the Sasaki type d--metric for
   N--anholonomic spaces satisfy the bulk of properties of the metric
   and connection on the Riemannian space but "slightly"
   nonholonomically modified). }
\begin{definition}\label{dncdg} 
A spectral d--triple satisfying the mentioned seven
  Conditions for his h-- and v--irreversible components is a real one
  which defines a (d--spinor) N--anholonomic noncommutative geometry
  defined  by the data 
$(\mathcal{A}, \mathcal{H}_{N},\ \widehat{\D},\ J,\  \Gamma ^{[cr]}\ )$ 
 and derived for the  Dirac d--operator (\ref{ddo}). 
\end{definition}
For a particular case, when the N--distinguished structures are of Lagrange
(Finsler) type, we can consider:
\begin{definition} \label{dnclg} A spectral d--triple satisfying the
  mentioned  seven
  Conditions for his h-- and v--irreversible components is a real one
  which defines a Lagrange, or Finsler, (spinor)
  noncommutative geometry  defined by the data
 $(\mathcal{A}, \mathcal{H}_{(SL)},\ _{(L)}\widehat{\D},\ J,\ \Gamma ^{[cr]}\ )$
 and derived for the  Dirac d--operator (\ref{dlo}). 
\end{definition} 
In Ref. \cite{vncg}, we used the concept of d--algebra
$\mathcal{A}_{N} \doteq (\ ^{-}\mathcal{A},\ ^\star \mathcal{A})$ which we
 introduced  as a "couple" of algebras for respective h-- and
 v--irreducible decomposition of constructions defined by the
 N--connection. This is possible if 
$\mathcal{A}_{N} \doteq \ ^{-}\mathcal{A}\oplus\  ^\star \mathcal{A}),$ but
we can consider arbitrary noncommutative associative algebras
$\mathcal{A}$ if  the splitting is defined by the Dirac d--operator.

\subsection{Distance in d--spinor spaces: Main Result 3}

We can select N--anholonomic and Lagrange commutative geometries from
the corresponding Definitions \ref{dncdg} and \ref{dnclg} if we put
respectively $\mathcal{A}\doteq C^\infty (\mathbf{V})$ and
 $\mathcal{A}\doteq C^\infty (V_{(L)})$ and consider real spectral
 d--triples. One holds:
\begin{theorem} \label{mr3} {\bf  (Main Result 3)} Let
 $(\mathcal{A}, \mathcal{H}_{N},\ \widehat{\D},\ J,\  \Gamma ^{[cr]}\
  )$ \\  
  (or $(\mathcal{A}, \mathcal{H}_{(SL)},\ _{(L)}\widehat{\D},\ J,\
  \Gamma ^{[cr]}\ )$) defines a noncommutative geometry being
  irreducible for $\mathcal{A}\doteq C^\infty (\mathbf{V})$ (or 
 $\mathcal{A}\doteq C^\infty (V_{(L)})),$ where $\mathbf{V}$ (or
  $V_{(L)}$) is a compact, connected and  oriented  manifold without
  boundaries, of   spectral dimension $dim\ \mathbf{V}=n+m$
 (or $dim\ V_{(L)} =n+n$ ). In this case, there are
  satisfied the conditions:
\begin{enumerate} 
\item There is a unique d--metric  $\mathbf{g}(\widehat{\D} ) = (g,\ ^\star g)$ 
 of type ((\ref{metr})) on  $\mathbf{V}$\ (or  of type (\ref{slm}) on
 $V_{(L)})$ with the "nonlinear" geodesic distance defined by 
\begin{equation} \label{fngdd}
 d(u_1,u_2) = \sup _{f\in C(\mathbf{V})} \left\{ f(u_1,u_2) / \parallel [\D ,
 f ]\parallel \leq 1 \right\}
\end{equation}
 (we have to consider $f\in C(V_{(L)})$ and $\ _{(L)}\widehat{\D}$ if we compute
 $d(u_1,u_2)$ for Lagrange  configurations). 
\item The N--anholonomic manifold $\mathbf{V}$\ (or Lagrange space
  $V_{(L)})$ is a spin N--anholonmic space (or a spin Lagrange
  manifold) for which the operators $\D ^\prime$ satisfying
  $\mathbf{g}(\D ^\prime) = \mathbf{g}(\widehat{\D})$ define an union of affine
  spaces identified by the d--spinor structures on $\mathbf{V}$ (we
  should consider the operators  $\ _{(L)}\D ^\prime$ satisfying
  $\ ^{(L)}\mathbf{g}(\ _{(L)}\D ^\prime) =\   ^{(L)}\mathbf{g}(\
  _{(L)}\widehat{\D})$ for the space   $V_{(L)})$).\
\item The functional $S(\D ) \doteq \int |\D |^{-n-m+2}$
    defines a   quadratic d--form with $(n+m)$--splitting
 for every affine spaces which is minimal for  
 $\widehat{\D} = \overleftarrow{\D}$ as the Dirac d--operator
  corresponding to the d--spin structure with the minimum proportional
  to the Einstein--Hilbert action constructed for the canonical
  d--connection with the d--scalar curvature $
  \overleftarrow{\mathbf{R}}$ (\ref{sdccurv}),
\footnote{The  integral for the  usual Dirac operator related to
  the Levi--Civita     connection  $D\ $ is computed:\ 
 $\int |D|^{-n+2} \doteq \frac{1}{2^{[n/2]}{\Omega _n}} Wres
    |D|^{-n+2},$ 
 where $\Omega _n$ is the integral of the  volume on the sphere
  $S^n$ and $Wres$ is the Wodzicki residu, see  details in Theorem 7.5
  \cite{bondia}. On N--anholonomic manifolds, we may consider
   similar definitions and computations but
    applying N--elongated partial derivatives and differentials. } 
\begin{equation*}
  S(\overleftarrow{\D}) = - \frac{n+m-2}{24}\  \int _{\mathbf{V}}\   
 \overleftarrow{\mathbf{R}}\ \sqrt{g}\ \sqrt{h}\ dx^1...dx^n\   
\delta y^{n+1}... \delta y^{n+k}.
\end{equation*}
\end{enumerate}
\end{theorem}
\begin{proof} In this work, we sketch only the idea and the key points of
 a such Proof.  The Theorem is a generalization for N--anholonomic
 spaces of a similar one, formulated in Ref. \cite{connes1}, with a
 detailed proof prestented in  \cite{bondia}, which seems to be a
 final  almost generally  accepted result. There are
 also alternative considerations,  with useful details, in
 Refs. \cite{rennie1,lord}.  
 For the Dirac d--operators,  we have to start with the
 Proposition \ref{pdohv} and then to repeat all constructions from 
 \cite{connes1,bondia}, both on h-- and v--subspaces, in N--adapted
 form.

 The  existence of a canonical d--connection structure which is metric
 compatible and  constructed from the coefficients of the 
 d--metric and  N--connecti\-on  structure is a crucial result allowing
 the formulation and proof of the Main Results 1-3 of this
 work. Roughly speaking, if the commutative Riemannian geometry can be
 extracted from a noncommutative geometry, we can also generate (in a similar, but
 tecnically more sophisticate form)  Finsler like geometries and
 generalizations. To do this, we have to consider the 
 corresponding parametrizations of the nonholonomic frame structure,
 off--diagonal metrics and deformations of the linear connection
 structure, all constructions being adapted to the N--connection
 splitting. If a fixed d--connection satisfies the metricity conditions, the
 resulting Lagrange--Finsler geometry belongs to a  class
 of nonholonomic Riemann--Cartan geometies, which (in their turns) are
   equivalents, related by nonholonomic maps, of  Riemannian spaces, see
 \cite{v1,v3}. However, it is not yet clear how
 to perform a such general proof for nonmetric d--connections (of
 Berwald or Chern type).  We shall present the technical details of such
 considerations in our further works.

Finally, we emphasize that for the Main Result 3 there is the
possibility to elaborate an alternative proof (like for the Main
Result 2) by verifying that the basic formulas proved for the
Riemannian geometry hold true on N--anholonomic manifolds by a
corresponding substitution of the N--elongated differentil and
partiale derivatives operators acting on canonical d--connections and
d--metrics. All such constructions are elaborated in N--adapted form
by preserving the respective h- and v--irreducible decompositions. $\Box$
\end{proof}

Finally, we can formulate three important conclusions:
\begin{conclusion} The formula (\ref{fngdd})  defines the
  distance in a manner as to be satisfied all necessary properties (finitenes,
  positivity conditions, ...) discussed in details in
  Ref. \cite{bondia}. It allows to generalize the constructions for
  discrete spaces with anisotropies and to consider anisotropic
  fluctuations of noncommutative geometries \cite{mart2,mart3} (of Finsler type, and more
  general ones, we omit such constructions in this work).
 For the nonholonomic configurations we have to work
  with canonical d--connection and d--metric structures.
\end{conclusion}
 Following the N--connection formalism originally elaborated 
  in the framework of  Finsler geometry, we may state:
\begin{conclusion} In the particular case of the canonical N--connection,
  d--con\-nec\-ti\-on and d--metrics defined by a regular Lagrangian, it is
  possible  a noncommutative geometrization of Lagrange mechanics
  related to  corresponding  classes of noncommutative Lagrange--Finsler geometry. 
\end{conclusion}
Such geometric methods have a number of applications in modern gravity:
\begin{conclusion} By anholonomic frame transforms, we can generate
  noncommutative Riemann--Cartan and Lagrange--Finsler spaces, in
  particular exact solutions of the Einstein equations with
  noncommutative variables \footnote{see
  examples in Refs. \cite{v0,vs1, v2,v3,vncgs}}, by considering
  N--anholonomic deformations of the Dirac operator.
\end{conclusion}

\subsection*{Acknowledgment}

The  work  summarizes the  results communicated in  a series of  Lectures and
Seminars:  

\begin{enumerate}

\item S. Vacaru, A Survey of (Non) Commutative Lagrange and Finsler
Geometry, lecture 2: Noncommutative Lagrange and Finsler
Geometry.  Inst. Sup. Tecnico, Dep. Math., Lisboa,
Portugal, May 19, 2004 (host:  P. Almeida).

 \item S. Vacaru, A Survey of (Non) Commutative Lagrange and Finsler
 Geometry, lecture 1:  Commutative Lagrange and Finsler Spaces and
 Spinors. Inst. Sup. Tecnico, Dep. Math., Lisboa,
 Portugal, May 19, 2004 (host:  P. Almeida).

 \item S. Vacaru, Geometric Models in Mechanics and Field Theory, 
 lecture at the Dep. Mathematics, University of Cantabria, Santander, 
Spain, March 16, 2004 (host: 
 F. Etayo).

 \item S. Vacaru, Noncommutative Symmetries Associated to the Lagrange and
 Hamilton Geometry,   seminar at the Instituto de Matematicas y
 Fisica Fundamental,  Consejo Superior de Investigationes, Ministerio 
de Ciencia y Tecnologia, Madrid, Spain, March 10, 2004 (host: M. de Leon).

 \item S. Vacaru, Commutative and Noncommutative Gauge Models,   lecture 
at the Department of Experimental Sciences, University of Hu\-el\-va,
Spain, March 5, 2004 (host:  M.E. Gomez).

\item S. Vacaru, Noncommutative Finsler Geometry, Gauge Fields and Gravity,
  seminar at the Dep. Theoretical Physics, University of Zaragoza, Spain,
 November 19, 2003 (host: L. Boya).

\end{enumerate}

The author is  grateful to the Organizes and hosts  for financial support and
collaboration. He also thanks  R. Picken for help  and  P. Martinetti 
  for  useful  discussions. 

%\newpage

\end{document}